\documentclass[reqno,a4paper]{amsart}
\usepackage{amsmath,amssymb,amsfonts,epsfig,float,color,subfigure,caption,textcomp}
\usepackage[colorlinks=true,breaklinks=true,linkcolor=lightblue,citecolor=lightblue]{hyperref}
\definecolor{lightblue}{rgb}{0.22,0.45,0.70}
\definecolor{cgray}{rgb}{0.9,0.9,0.9}
\graphicspath{figures/}
\newcommand\vdiv{\mathop{\mathrm{div}}\nolimits}

\newcommand{\abs}[1]{\ensuremath{\left|#1\right|}}

\def\x{{\,\mathrm{x}\x}}

\def\L{\mathrm{L}}

\numberwithin{equation}{section}

\renewcommand\L{\mathrm{L}}

\renewcommand\L{\mathrm{L}}

\newtheorem{theorem}{Theorem}[section]
\newtheorem{lemma}[theorem]{Lemma}

\newtheorem{proposition}{Proposition}[section]
\def\exp{\mbox{{\rm exp}}}

\def\#{{{\cal D}_h}}
\usepackage{bbold}
\usepackage{tikz}
\usepackage{pgfplots}
\usetikzlibrary{spy,calc}

\title[Kinetic derivation of a time-dependent SEIRD reaction-diffusion system for COVID-19]{Kinetic derivation of a time-dependent SEIRD reaction-diffusion system for COVID-19 }

\author{\textbf{Mohamed Zagour}}

\address{Institut National des Sciences Appliqu\'{e}es, Universit\'{e} EuroMed F\`{e}s, Maroc}
\email{zagourmohamed@gmail.com}
\begin{document}
\maketitle
\begin{abstract} In this paper, we propose a time-dependent Susceptible-Exposed-Infectious-Recovered-Died (SEIRD) reaction-diffusion system for the COVID-19 pandemic and we deal with its derivation from a kinetic model. The derivation is obtained by mathematical description delivered at the micro-scale of individuals. Our approach is based on the micro-macro decomposition which leads to an equivalent formulation of the kinetic model which couples the microscopic equations with the macroscopic equations. We develop a numerical asymptotic preservation scheme to solve the kinetic model. The proposed approach is validated by various numerical tests where particular attention is paid to the Moroccan situation against the actual pandemic.
	
\end{abstract}
%%%%%%%%%%%%%%%%%%%%%%%%%%%%%%%%%%%%%%%%%%%%%%%%%%%%%%%%%%%%%%%%%%%%%%%%%%%%%%%%%%%%%%%%%%%%%%%%%%%%%%%%%%%%%
\section{Introduction}
In December 2019, the outbreak of the new coronavirus called COVID-19 caused by severe acute respiratory syndrome coronavirus 2 (SARS-CoV-2) first occurred in Wuhan, China. It is spreading so fast and has spread to almost other countries around the world. Therefore, the World Health Organization declared it to be a pandemic since March 11, 2020. Unfortunately, this pandemic so severely affects the economy, health, security of society. Currently, as of January 17, 2021 more than 95,013,471 people are infected and more than 2,032,253 have died from this virus, see update data in \cite {[Rf1]}.

As it is known, mathematical models can help in many tasks, for example verifying different epidemic scenarios, estimating transmission parameters, testing various hypotheses, and better understanding the mechanisms of contagion. In fact, it provides decisions to optimize possible control strategies, such as containment measures, lockdowns, and vaccination campaigns. However, several models have been proposed to describe the evolution of epidemics which can be collective models or network models. We mention that collective models describe the spread of the epidemic in a population using a limited number of collective variables with a small number of parameters. We can find logistic models \cite{[KM27], [V38]}, generalized growth models \cite{[C17]}, Richards models \cite{[R59]}, Susceptible-Infected-Recovered models (SIR) \cite{[B75], [KM27]} and SEIR (Susceptible-Exposed-Infectious-Removed) models \cite{[C17]}. Note that SIR, SEIR, and other similar models belong to the class of compartmental models, see \cite{[B17],[C17],[R08]}. On the other hand, network models treat a population as a network of interacting individuals, and the contagion process is described at the microscopic level, see for example \cite{[K05], [W14]}.

 Recently, several mathematical models have been proposed and developed to understand the dynamics of COVID-19, see e.g \cite{[AI20],[BB20],[CSB20],[Est2020],[FP20],[Fox20],[LT20],[LZ20],[MK20],[PLZ20],[SZ20],[SW20],[ZLY20]}. However, in this paper we are interested in a collective modified SEIRD model. Specifically, the improved model  is given by the following system in nondimensional form
\begin{equation}\label{Cross-Diffusion}
\left\{\begin{array}{l}
\displaystyle\partial_t S -d_1\Delta S = A-\mu S-\beta(t)S\frac{I}{N}, \\ \\
\displaystyle\partial_t E -d_2\Delta E = \beta(t)S\frac{I}{N}-(\mu+\xi)E, \\ \\
\displaystyle\partial_t I -d_3\Delta I = \xi E-(\gamma+\mu+\alpha)I, \\ \\
\displaystyle\partial_t R -d_4\Delta R = \gamma I-\mu R, \\ \\
\displaystyle\partial_t D  = \alpha I,
\end{array}
\right.
\end{equation}
in $\Omega_T :=(0,T)\times \Omega$ for a fixed time $T>0$.
We augment this system along with the boundary conditions 
\begin{equation}
d_{1} \nabla S\cdot\;\eta=d_{2} \nabla E\cdot\;\eta=d_{2} \nabla I\cdot\;\eta=d_{3} \nabla R\cdot\;\eta=0\quad  \mbox{on }\;  \Sigma_T:=(0,T]\times\partial\Omega,
\end{equation}
and the initial conditions
\begin{equation*}
S(t=0,x)=S_{0}(x),\; E(t=0,x)=E_{0}(x),\;I(t=0,x)=I_{0}(x),
\end{equation*}
\begin{equation}
R(t=0,x)=R_{0}(x),\;D(t=0,x)=D_{0}(x),\; \quad \mbox{for} \;\; x \in \Omega,
\end{equation}

\begin{figure}[h!]
	\centering
	\subfigure{\includegraphics[height=1.7in ,width=4.5in]{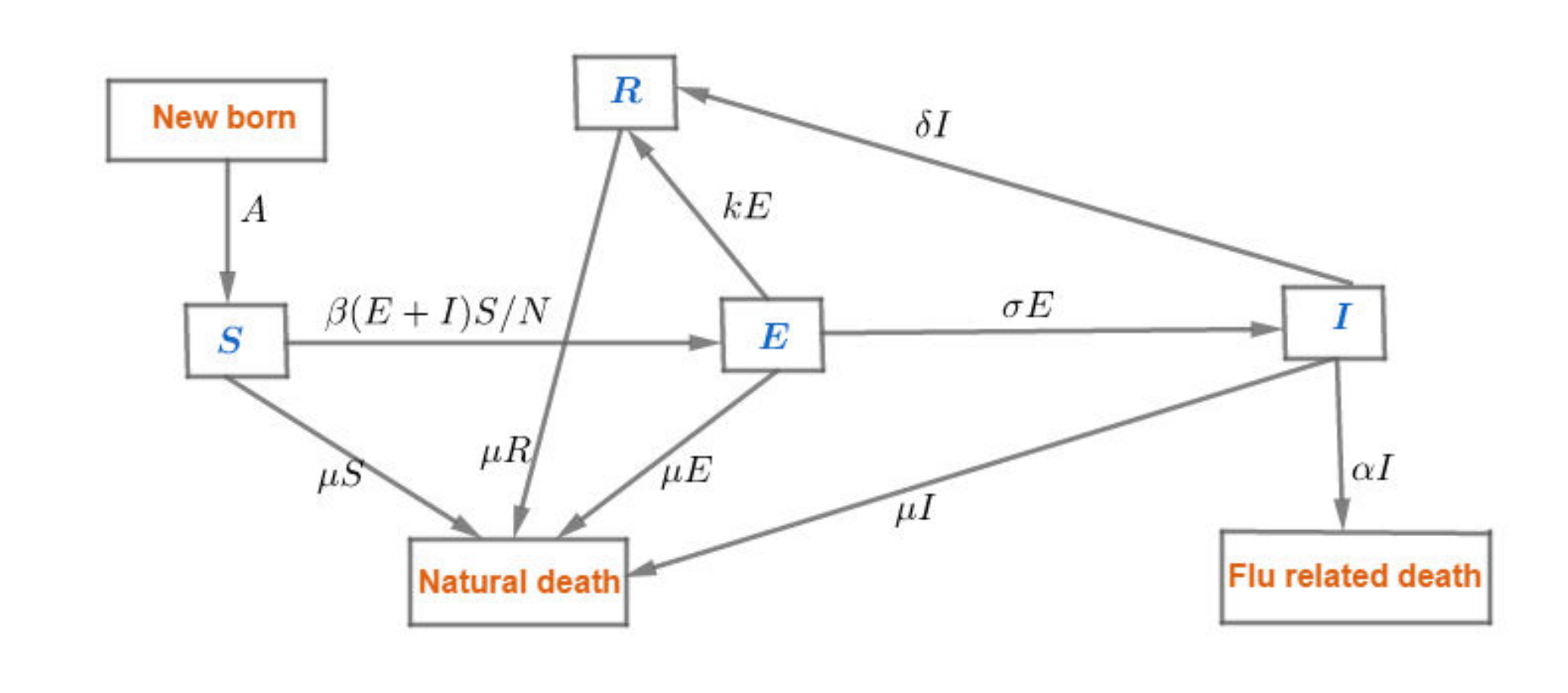}}
	\caption{The flow diagram of spatially homogeneous SIERD model}
	\label{Flow}
\end{figure}

In all of the above equations, the constants $ d_i; \, i = 1, \cdots, 4 $ are the diffusivity constants. A description of all variables and parameters used in the aforementioned system \eqref{Cross-Diffusion} is presented in tables \ref {tab1} and \ref {tab2} respectively. The scheme is illustrated in Figure \ref{Flow}. Note that if $ A = \mu = 0 $ and $ \ xi = \infty $, the \eqref{Cross-Diffusion} system is reduced to the classic SIR system which means that it has no latent period. Let us now briefly describe the significance of the variables: $ S $ is the number of individuals likely to be exposed, while $ E $ is the number of individuals exposed, this number constitutes a class where the disease is latent and the individuals are infected and non-infectious. Different processes can occur. For example, Susceptible may become Exposed due to contact with infectious individuals and transmission function. Additionally, Exposed can become infectious with a $ \xi $ rate and infectious recover with a $ \gamma $ rate.
We recall that we consider a time-dependent transmission rate function, which incorporates the impact of government action (i.e total or partial lockdown), wearing the mask and respecting sanitary protocol such as the social distance. A typical example is a step-wise function, see Section \ref{Sec4}.
%\begin{equation}
%\beta(t)=\beta_0\big((1-\phi)e^{-qt}+\phi\big),
%\end{equation}
%Thus, this time-dependent transmission rate $\beta(t)$ allows the model to capture early sub-exponential growth dynamics whenever $R_0 > 1$ and $q > 0$. A special case is the standard SIR
%modeling framework when $q = 0$. The time-dependent transmission rate $\beta(t)$ declines exponentially from an initial value $\beta_0$ towards $\phi\beta_0$ at the rate $q>0$. Assuming that $R_0 > 1$ and the population size $N$ is sufficiently large such that the effect of susceptible depletion on the effective reproduction number due to new infections is negligible, the quantity $(1-\phi)$ models the proportionate reduction in $\beta$ that is needed to reach an effective stationary reproduction number of 1. Hence, $\phi$ can be simply estimated as $\frac{\gamma}{\beta_0}$.
\begin{table}[h!]
	\begin{center}
		\caption{Description of the variables of the SEIRD system}
		\label{tab1}
		\begin{tabular}{l|l} % <-- Changed to S here.
			\textbf{Variable } & \textbf{Description} \\
			\hline
			$N(t,x)$ & Total number of live individuals \\
            %\hline
            $S(t,x)$ & Population of susceptible individuals \\
            %\hline
            $E(t,x)$ & Population of exposed individuals \\
            %\hline
            $I(t,x)$ & Population of infected individuals \\
            %\hline
            $R(t,x)$ & Population of  recovered individuals \\
            %\hline
            $D(t,x)$ & Population of  died individuals \\
            \hline
		\end{tabular}
	\end{center}
\end{table}

\begin{table}[h!]
	\begin{center}
		\caption{Description of the parameters of the SEIRD system}
		\label{tab2}
		\begin{tabular}{l|l} % <-- Changed to S here.
			\textbf{Parameter} & \textbf{Description} \\
			\hline
			$A$ & Recruitment rate assumed $A=\mu\,N$\\
			%\hline
			$\mu$ & Natural death rate for susceptible individuals\\
			%\hline
%			$\beta_0$ & Initial value of transmission rate\\
			%\hline
			$\xi$ &Rate of progression from exposed to infectious\\
			%\hline
			$\gamma$ & Recovery rate of infectious individuals  \\
			%\hline
			$\alpha$ &Virus-induced average fatality rate \\
			\hline
		\end{tabular}
	\end{center}
\end{table}

%===========================Rf makro===============

Note that in the absence of the diffusion terms ($d_i=0,\; i=1,\cdots 4$) and that the transmission rate is constant (i.e $\beta(t)=\beta$), system \eqref{Cross-Diffusion} reduces to a spatially homogeneous SEIRD problem, see for e.g. \cite{[AT04],[DHB13],[H00],[HCZ20],[KR11],[PNW98]}.  However, if we assume that $d_i\neq0$ we obtain a spatially inhomogeneous SEIRD problem, among others we cite \cite{[Mass07],[SSL10]}. Comparing to the mentioned papers above, in this paper we propose a time-dependent SEIRD reaction-diffusion system, where the aim is to take into account the fast change in time of the transmission rate function. In fact, this help to well model the different strategies taken to defeat the virus, for instance partial or total lockdown and now the vaccination. In passing, we would like to mention that the basic reproduction ratio, denoted by $R_0$, is the classical epidemiological measure associated with the reproductive power of the disease. It is used to estimate the growth of the viral epidemic. For our system \eqref{Cross-Diffusion} it is given by
\begin{equation}
R_0(t)=\frac{\xi}{(\xi+\mu)(\gamma+\alpha+\mu)}\,\beta(t).
\end{equation}
Note that $R_0(t)$ provides a threshold for disease-free equilibrium point stability. Precisely, if $R_0(t) < 1$, the disease goes out; while if $R_0(t)> 1$, an epidemic occurs, see e.g. \cite{[H00],[VW02]}.

In this paper, we deal with the derivation of the time-dependent SEIRD reaction-diffusion system \eqref{Cross-Diffusion} from kinetic model by using the micro-macro decomposition method. The main idea of this method consists in rewriting the kinetic model as coupled system of microscopic part and macroscopic one.  Many works used this method within different fields of application. For instance, chemotaxis phenomena in the basis of the famous Keller-Segel model \cite{[BBC16]}, formation of patterns induced by cross-diffusion in a fluid \cite{[ABKMZ19],[BKZ18]}. In fact, this technique has been adopted to design a numerical scheme that preserves the asymptotic property introduced by \cite{[JI99],[KL98]}. In other words, a stable numerical scheme in the limit along the transition from kinetic to macroscopic regimes. 

This paper is organized as follows. In Section \ref{Sec2} we summarize the micro-macro method which leads to the derivation of system \eqref{Cross-Diffusion} from a kinetic model. Section \ref{Sec3} is devoted to the development of an asymptotic preserving numerical (AP)-scheme in one dimension, inspiring from the paper by \cite{[ABKMZ19]}. In other words, the uniform stability with respect to the parameter $\varepsilon$ and the consistency with the reaction-diffusion limit. In Section \ref{Sec4}, we present the numerical simulations obtained from micro-macro formulation and from the macroscopic scheme, where we show the asymptotic preserving scheme property. Moreover, we demonstrate the effect of presence of the diffusion terms in system \eqref{Cross-Diffusion}, and its sensitivity with respect to the different choices of the reproduction ratio $R_0$. Finally, particular attention is paid to the Moroccan situation against the actual pandemic.

%%%%%%%%%%%%%%%%%%%%%%%%%%%%%%%%%%%%%%%%%%%%%%%%%%%%%%%%%%%%%%%%%%%%
\section{From kinetic model to  SEIRD reaction-diffusion system}\label{Sec2}
%%%%%%%%%%%%%%%%%%%%%%%%%%%%%%%%%%%%%%%%%%%%%%%%%%%%%%%%%%%%%%%%%%%
This section aims to derive briefly SEIRD reaction-diffusion system \eqref{Cross-Diffusion} from  kinetic model using micro-macro decomposition method by following the line of paper \cite{[ABKMZ19]}. We present the properties of the aforesaid model. On the basis of the micro-macro decomposition technique, we give an equivalent appropriate system. 
\subsection{Kinetic model}
This subsection is devoted to state kinetic model and to present its properties. We consider the following  kinetic model for $i=1,\dots,4$
\begin{equation}\label{so}
\left\{
\begin{array}{l}
\displaystyle\varepsilon \partial_tf_i+ v  \cdot \nabla_{x} f_i=
\frac{1}{\varepsilon}\mathcal{T}_i(f_i) +G_i(f_1,\dots       ,f_4),    \\\\
  \displaystyle\partial_t D=\alpha \int_Vf_3\,dv,\\\\
f_i(t=0,x,v)=f_{i,0}(x,v),\qquad D(t=0,x)=D_0,
\end{array}\right.
\end{equation}
where $f_1(t,x,v),\cdots, f_4(t,x,v)$ are the distribution functions describing the statistical evolution of susceptible, exposed, infected and recovered individuals, respectively. $t>0$, $ x\in \mathbb{R}^{d}$, $v \in V$ are respectively, time, position and velocity. The term  $\mathcal{T}_i$ is the stochastic operator representing a random modification of direction of individuals and the operator $G_i$ ($i=1,\dots ,4$) describing the gain-loss balance of individuals. The mean free path $\varepsilon$ measures the distance between individuals. We mention that we adopt the parabolic-parabolic scaling limit, see for more details \cite{[ABKMZ19]}. \\
The micro-macro decomposition method is based on the following assumptions. The turning operator $\mathcal{T}_i$ are given by
\begin{equation}
\displaystyle\mathcal{T}_{i}(f_i)= \int_{V}\big(T^{*}_{i}(v^*,v)f_i(t, x, v ^{*}) -
T_{i}(v,v^*) f_i(t, x, v ) \big)dv ^{*}, \label{3.4}
\end{equation}
where $T_i$ is the probability kernel for the new velocity $v \in V$ given that the previous velocity was $v ^*$. We assume that the operators $\mathcal{T}_i$ satisfy
\begin{equation}\label{H0}
\displaystyle \int_V \mathcal{T}_i\,dv  =0, \;\;i=1,\dots ,4.
\end{equation}
Moreover, we assume that there exists a bounded velocity distribution $M_i(v)>0$ for $i=1,\dots      ,4$ independent of
$t$ and $x$ such that
\begin{equation}\label{tx}
 T_i (v ,v ^{*} ) M_i(v ^{*}) =  T_i (v ^{*},v  ) M_i(v),
\end{equation}
holds. The flow produced by these equilibrium
distributions vanishes and $M_i$ are normalized, i.e.
\begin{equation}\int_V v  \, M_i(v )dv   =0, \quad \int_V
M_i(v ) dv  =1, \quad  i=1,\dots      ,4. \label{equilibre}\end{equation}\\
Regarding the probability kernels, we assume that $T_i^1(v ,v ^{*})$ is bounded,
and there exist a constant $\sigma_i >0$ ($ i\in\{1,\dots      ,4\}$),
such that
\begin{eqnarray}\label{cx}
T_i^1(v ,v ^{*})\geq \sigma_i M_i(v ),
\end{eqnarray}
for all $ (v ,v ^{*}) \in V\times V $,  $ x \in \Omega $ and $t>0$.
Next, the interaction operators $G_i$ satisfy the following properties	
\begin{equation}\label{CO}\displaystyle\int_V G_{i}(f_1,\dots       ,f_4,v) dv  =0,\;\;i = 1,\dots,4.
\end{equation}

\noindent Using the same arguments as in \cite{[ABKMZ19]}, the operator $\mathcal{T}_i$ has the following properties.
\begin{lemma}
	\label{LE1}  By assuming that the hypothesis \eqref{H0}, \eqref{tx} and \eqref{equilibre} are satisfied. Then, the
	following properties of the operators $\mathcal{T}_i$ for $i=1,\dots,4$ hold true:
	\begin{itemize}\label{L1}
		\item[i)] The operator $\mathcal{T}_i$ is self-adjoint in the space
		$\displaystyle{{\L^{2}\left(V ,{dv \over M_i(v)}\right)}}$. 
		\item[ii)]
		For $f\in \L^2$, the equation $\mathcal{T}_i(g) =f$ has a unique
		solution $\displaystyle{g \in \L^{2}\left(V,\frac{dv}{
				M_i(v)}\right)}$, satisfying
		$$\int_Vg(v) dv  = 0 \quad \Longleftrightarrow \quad   \int_{V} f(v )\, dv  =0. 
		$$
		\item[iii)]  The equation $\mathcal{T}_i(g) =v  \,  M_i(v)$, has a
		unique solution denoted by $\theta_i(v )$ for $ i=1,\dots      ,4$. 
		\item[iv)]  The kernel
		of $\mathcal{T}_i$ is $N(\mathcal{T}_i) = vect(M_i(v))$ for $ i=1,\dots      ,4$.
	\end{itemize}
\end{lemma}
%%%%%%%%%%%%%%%%%%%%%%%%%%%%%%%%%%%%%%
\subsection{Micro-macro formulation}
%%%%%%%%%%%%%%%%%%%%%%%%%%%%%%%%%%%%%%
This subsection is devoted to obtain an equivalent micro-macro system of nonlocal kinetic-fluid model (\ref{so}). The obtained equivalent system contains microscopic and macroscopic components. 
%In what follows, the integral with respect to the variable $v $ will be denoted by $\langle . \rangle$.

\noindent The main idea of the micro-macro method is to decompose the distribution function $f_i$ for $i=1,\dots      ,4$ as follows
$$f_i(t,x,v )=M_i(v) u_i(t,x) + \varepsilon  g_i(t,x,v), $$
where 
$$u_i(t,x)= \langle f_i(t,x,v )\rangle:=\int_Vf_i(t,x,v )\,dv .$$ 
We have $\langle g_i \rangle=  0$ for $ i=1,\dots      ,4$. Inserting $f_i$ in kinetic model (\ref{so}) and using the above assumptions and properties of the interaction and the turning operators, we get
\begin{equation}\label{n} 
\displaystyle\partial_t (M_i (v)u_i)  + \varepsilon \partial_t g_i +
\frac{1}{\varepsilon}  v  M_i(v ) \cdot \nabla u_i + v  \cdot
\nabla g_i   = \frac{1}{\varepsilon}\mathcal{T}_i(g_i)+G_{i}(f_1,\dots,f_4 ),
\end{equation}

\noindent The micro-macro decomposition method is based on two steps. First, we use the projection technique to separate the macroscopic density $u_i(t,x)$ and microscopic quantity $g_i(t,x,v )$ for $i=1,\dots      ,4$. For that, we consider $P_{M_i}$ the orthogonal projection onto $N(\mathcal{T}_i)$, for $i=1,\dots      ,4$. It follows
$$P_{M_i(v)}(h)= \langle h\rangle M_i(v), \quad  \mbox{for any}\quad  h\in
\displaystyle{{\L^{2}\left(V ,{ dv  \over M_i(v )}\right)}}, \qquad i=1,\dots      ,4.$$ 

%\noindent Regarding the orthogonal projections $P_{M_i(v)}$ for $i=1,\dots      ,4$, we have the following result.

%\begin{lemma}  \label{LE2} We have the following properties for the projection $P_{M_i(v )}$, $ i=1,\dots      ,4$
%	$$ (I-P_{M_i(v)})\big(M_i(v) u_i\big)= P_{M_i(v )}(g_i)=0, $$
%	$$(I-P_{M_i(v )})\big(v  M_i(v) \cdot\nabla u_i) u_i\big)=v  M_i(v )\cdot \nabla u_i,$$	
%	$$(I-P_{M_i(v)})(\mathcal{T}_i(g_i))=\mathcal{T}_i(g_i)$$	
%\end{lemma}
Now, inserting the operators $I-P_{M_i}$ into Eq.  (\ref{n}), and using known properties for the projection
$P_{M_i}\;i= 1,\cdots4$ yields the following micro-macro formulation
\begin{equation}\label{mM1}
\left\{
\begin{array}{l l}
\displaystyle\partial_t g_i +
\frac{1}{\varepsilon^2} v  M_i(v) \cdot \nabla u_i+ \frac{1}{\varepsilon}(I-P_{M_i})(v 
\cdot
\nabla g_i) =\frac{1}{\varepsilon^2}\mathcal{T}_i(g_i)+ \frac{1}{\varepsilon}(I-P_{M_i}) G_{i}(f_1,\dots  ,f_4), \\
{}\\
\displaystyle\partial_t u_i+ \langle v    \cdot
\nabla g_i \rangle =  \langle G_{i}(f_1,\dots,f_4)\rangle, \\\\
\displaystyle\partial_t D=\alpha u_3,
\end{array} \right.
\end{equation}

The following proposition shows that micro-macro formulation \eqref{mM1} is equivalent to nonlocal kinetic-fluid equation (\ref{so}).

\begin{proposition}	
	\noindent i) Let $(f_1,\dots    ,f_4)$ be a solution of kinetic model (\ref{so}). Then \\$(u_1,\dots,u_n,g_1,\dots    ,g_4)$ (where
	$u_i=\langle f_i \rangle $ and $g_i=
	{1\over \varepsilon}(f_i-M_i u_i)$) is a solution to coupled system
	(\ref{mM1}) associated with the following initial data for $i=1,\dots    ,4$
	\begin{equation} \label{er4} u_i(t=0)=u_{i,0} =\langle f_{i,0} \rangle, \quad
	g_i(t=0)=g_{i,0}={1\over \varepsilon}(f_{i,0}-M_i u_{i,0}) 
	\end{equation}	
	\noindent ii) Conversely, if $(u_1,\dots    ,u_n,g_1,\dots    ,g_4)$ satisfies system
	(\ref{mM1}) associated with the following initial data\\ $(u_{1,0},\dots    ,u_{n,0}, g_{1,0},\dots    ,g_{n,0})$ such
	that $\langle g_{i,0} \rangle=0$ for $i=1,\dots    ,4$. Then $(f_1,\dots    ,f_4)$ (where $f_i=M_i u_i+\varepsilon g_i$)  is a solution to  kinetic model (\ref{so}) with initial data
	$f_{i,0}=M_i u_{i,0}+\varepsilon g_{i,0}$ and we have $ u_i=\langle f_i \rangle$ and $\langle g_i\rangle=0$, for $i=1,\dots    ,4$.
\end{proposition} 

\noindent Next, in order to develop asymptotic analysis of
system (\ref{mM1}), we assume that $G_{i}$ satisfy the following asymptotic behavior $\varepsilon \to 0$
\begin{equation} \label{G} G_{i}\Big(M_1(v )u_1 +\varepsilon g_1,\dots    , M_4(v )u_4 +\varepsilon g_4\Big)= G_{i}^{j}\Big(M_1(v)u_1,\dots    , M_4(v )u_4 \Big)+ O(\varepsilon), \end{equation}
for $ i=1,\dots    ,4$ .
Now, we show that micro-macro formulation (\ref{mM1}), which is equivalent to kinetic equation (\ref{so}), allows to obtain a general macroscopic model as $\varepsilon$ goes to $0$. Indeed, using  (\ref{G}) and (\ref{mM1}), we obtain for $i=1,\dots    ,4$
$$\mathcal{T}_i(g_i) =  v  M_i(v) \cdot \nabla  u_i +O(\varepsilon).$$

\noindent From Lemma \ref{L1}, property $ii)$, the operator $\mathcal{T}_i$ is invertible. This implies 
\begin{eqnarray}\label{x1} g_i =  \mathcal{T}_i^{-1}\Big(v  M_i \cdot \nabla u_i\Big)+O(\varepsilon),\;\;i=1,\dots    ,4.
\end{eqnarray}
Next, inserting (\ref{x1}) into the second equation in (\ref{mM1}) yields the following macroscopic system
\begin{equation}\label{G1}
\displaystyle \partial_t u_i +  \Big\langle v  \cdot \nabla \mathcal{L}_i^{-1}\Big(v  M_i(v) \cdot \nabla  u_i\Big)\Big\rangle 
= \Big\langle G_{i}(M_1(v ) u_1,\dots    , M_4(v )u_4)\Big\rangle + O(\varepsilon),
\end{equation}
We have the following 
$$\displaystyle\left\langle v \cdot \nabla
\mathcal{T}_i^{-1}\Big(v  M_i(v ) \cdot \nabla  u_i\Big)\right\rangle= \vdiv\Big( \left\langle v  \otimes
\theta_i(v )\right\rangle\cdot \nabla u_i \Big), $$ 
where
	$\theta_i(v)$ are given in Lemma \ref{LE1} for $ i=1,\dots    ,4$.

\noindent Finally, we obtain the following reaction-diffusion system
\begin{equation}\label{mM2}
\left\{
\begin{array}{l l}
\partial_t u_i -  \vdiv \, \Big(D_i \cdot \nabla u_i\Big)=H_i(u_1,\dots    ,u_4) +O(\varepsilon),\
\\\\
\displaystyle\partial_t D=\alpha u_3,
\end{array}
\right.
\end{equation}
\noindent where the functions $D_i$ and $H_i$ are given by
\begin{equation}\label{di}
\qquad D_i=- \int_Vv  \otimes \theta_i(v)  dv  , 
\end{equation}

\begin{equation} \label{H}H_i(u_1,\dots,u_4)= \int_VG_{i}(M_1 u_1,\dots    , M_4 u_4) dv  ,\;  \hbox{for}\; i=1,\dots    ,4. \end{equation}
Now, we consider in (\ref{so}) a particular choice of terms:
$$u_1=S, \quad u_2=E,\quad u_3=I,\quad u_4=R. $$
Next, we assume that the probability kernel $T_i$ is given by 
\begin{equation*}
T_i= \frac{\sigma_i}{M_i(v)},\quad \hbox{for} \; i=1,\dots    ,4.
\end{equation*} 
This implies  
\begin{equation}\label{OpL}
\mathcal{T}_i (g)= -\sigma_i \Big( g -M_i \langle g \rangle \Big)= -\sigma_i\;g\quad \hbox{for} \; i=1,\dots    ,4.\end{equation} 
Using (\ref{equilibre}), (\ref{OpL}) and Lemma \ref{LE1}, then $\theta_i$ is given by
$$\theta_i= - \frac{1}{\sigma_i}v  M_i(v).$$

\noindent The modeling of the interaction operators $G_{i}$ is given by 
\begin{equation}\label{GG1}
\left\{
\begin{array}{l}
\displaystyle
\displaystyle G_{1}(f_1,\dots    ,f_4)=  \frac{1}{\abs{V}}\big(A-\mu f_1-\beta(t) f_1f_3/n\big),\\
\displaystyle G_{2}(f_1,\dots    ,f_4)=  \frac{1}{\abs{V}}\big(\beta(t) f_1f_3/n-(\mu+\xi)f_2\big),\\
\displaystyle G_{3}(f_1,\dots    ,f_4)=  \frac{1}{\abs{V}}\big(\xi f_3-(\gamma+\mu+\alpha)f_3\big),\\
\displaystyle G_{4}(f_1,\dots    ,f_4)=  \frac{1}{\abs{V}}\big(\gamma f_3-\mu f_4\big),
\end{array}
\right.
\end{equation}
where $n=\displaystyle\sum_{i=1}^{4}f_i$.
Then, we use the definition of $H_i$ in (\ref{H}) to obtain from (\ref{GG1})
\begin{equation} \label{HH}H_i(u_1,\dots    ,u_n)= F_i(u_1,\dots    ,u_n).\end{equation}
Finally, collecting the previous results and (\ref{mM2}), we obtain SEIRD  system \eqref{Cross-Diffusion} with diffusion of the order $O(\varepsilon)$
\begin{equation}\label{Cross-Diffusion2}
\left\{\begin{array}{l}
\displaystyle\partial_t S -d_1\Delta S = A-\mu S-\beta(t)S\frac{I}{N}+O(\varepsilon), \\ \\
\displaystyle\partial_t E -d_2\Delta E = \beta(t)S\frac{I}{N}-(\mu+\xi)E+O(\varepsilon), \\ \\
\displaystyle\partial_t I -d_3\Delta I = \xi E-(\gamma+\mu+\alpha)I+O(\varepsilon), \\ \\
\displaystyle\partial_t R -d_4\Delta R = \gamma I-\mu R+O(\varepsilon), \\ \\
\displaystyle\partial_t D  = \alpha I,
\end{array}
\right.
\end{equation}
%%%%%%%%%%%%%%%%%%%%%%%%%%%%%%%%%%%%%%%%%%%%%%%%%%%%%%%%
\section{Numerical methods}\label{Sec3}
%%%%%%%%%%%%%%%%%%%%%%%%%%%%%%%%%%%%%%%%%%%%%%%%%%%%%%%%
In this section we develop an asymptotic preserving numerical schemes(AP) in one dimension. In other words, the uniform stability with respect to the parameter $\varepsilon$ and the consistency with the reaction-diffusion limit. The discretization of problem (\ref{mM1})) is carried out with respect to each independent variable (time, velocity and space).
\subsection{Semi-implicit time discretization}
In this first step, we present a time discretization of our coupled system. We denote by $\Delta t$ a fixed time step, and by $t_k$  a discrete time such that $t_k=k\, \Delta t $  $ k\in N.$ The approximation of $u_1(t,x)$ and $g_i(t,x,v)$ at the time step $t_k$ are denoted respectively by $u_i^k\approx u_i(t_k,x)$ and $g_i^k\approx g_i(t_k,x,v)$.\\
In the first microscopic equation of (\ref{mM1}),  the only term which presents a stiffness in the collision part, for small $\varepsilon$, is $\displaystyle  \frac{1}{\varepsilon}\mathcal{T}_i(g_i)$. Hence,  we take an  implicit scheme to ensure the stability for this term, while  the other terms are still  explicit, then on has
\small\begin{equation}\label{mMD1}
\begin{array}{ll}
\displaystyle  \frac{ g_i^{k+1}-g_i^{k}}{\Delta t} +
\frac{1}{\varepsilon^2} v M_i \cdot \nabla u_i^k + \frac{1}{\varepsilon} (I-P_{M_i})(v\cdot\nabla g_i^k) =  \frac{1}{\varepsilon^{2}}\mathcal{T}_i(g_i^{k+1}) 
+\frac{1}{\varepsilon} (I-P_{M_i})G_i(u_1^k,u_2^k,u_3^k,u_4^k).
\end{array}
\end{equation}
In the second  macroscopic equation of (\ref{mM1}), we take $h$ at the time $t_{k+1}$, which gives
\begin{equation}\label{mMD2}
\frac{ u_i^{k+1}-u_i^{k}}{\Delta t}  + \langle v    \cdot
\nabla g_i^{k+1} \rangle = \left\langle G_i(u_1^k ,u_2^k,u_3^k,u_4^k)\right\rangle.
\end{equation}

\begin{proposition}
	The time discretization (\ref{mMD1})-(\ref{mMD2}) is consistent with (\ref{G1}) when $\varepsilon$  goes to $0$.
\end{proposition}

%%%%%%%%%%%%%%%%%%%%%%%%%%%%%%%%%%%%%%%%%
\subsection{Fully discrete scheme 1D}
In this section, we construct a suitable space discretization of (\ref{mMD1})- (\ref{mMD2}). The domain under consideration $[-L,L]$. The velocity space $[-V,V]$ can be treated by using a standard discretization. \\
We define a straggered grid $x_j=j \, \Delta x, $  $j=\,0,\, .\,.\,.,\, N_x$ with $\displaystyle N_x=\frac{L}{\Delta x}$, and the cell center points $x_{j-\frac{1}{2}}=(j-\frac{1}{2}) \Delta x, $  $j=\,0,\, .\,.\,.,\, N_x+1$. Let $u^k_{i,j} $ and $\displaystyle g^k_{i,j-\frac{1}{2}}$ be approximations of $u_i(t_k,x_j) $ and $\displaystyle g_i(t_k,x_{j-\frac{1}{2}}, v)$ respectively. Proceeding as \cite{[ABKMZ19],[BT18]}, the microscopic equation  (\ref{mMD1}) is discretized at points $\displaystyle x_{j+\frac{1}{2}}$ while the other macroscopic equations (\ref{mMD2})) are discretized at point $x_j$. Then, we obtain
\begin{equation}\label{mMDD1}
\begin{array}{ll}
\displaystyle  \frac{ g^{k+1}_{i,j+\frac{1}{2}}-g^{k}_{i,j+\frac{1}{2}}}{\Delta t} +
\frac{1}{\varepsilon^2}  v M(v)\frac{u^k_{i,j+1} -u^k_{i,j}}{\Delta x} + \frac{1}{\varepsilon} (I-P_{M_i})(v
\cdot \nabla_{x} g^k_{i,j+\frac{1}{2}} ) =  \frac{1}{\varepsilon^{2}}\mathcal{T}_i(h^{k+1}_{i,j+\frac{1}{2}}) \\\\
\hskip5.9cm+\frac{1}{\varepsilon}(I-P_{M_i})G_i(u^k_{1,j+\frac{1}{2}},\cdots,u^k_{4,j+\frac{1}{2}}),
\end{array}
\end{equation}
\begin{equation}\label{mMDD2}
\displaystyle
\frac{ u^{k+1}_{i,j}-u^{k}_{i,j}}{\Delta t}  + \langle v    \cdot
\nabla g^{k+1}_{i,j+\frac{1}{2}} \rangle = \langle G_i(u^k_{1,j},\cdots,u^k_{4,j})\rangle,
\end{equation}

\begin{proposition}\label{Pro3.2}
	The time and space approximation  (\ref{mMDD1} )-(\ref{mMDD2} ) of the  kinetic equation (\ref{so}) in the limit $\varepsilon$ goes to zero satisfy the following discretization
	\begin{equation}\label{mMD221}
	\begin{array}{ll}
	\displaystyle
	\frac{ u_i^{k+1}-u_i^{k}}{\Delta t}  + \frac{1}{\Delta x}\Big\langle  v\cdot  \Big [ \mathcal{T}_i^{-1}\Big( M(v) \cdot  \frac{u^k_{i,j+1} -u^k_{i,j} }{\Delta x} \Big) - \mathcal{T}_i^{-1}\Big(v M(v) \cdot  \frac{u^k_{i,j} -u^k_{i,j-1} }{\Delta x} \Big)\Big ] \Big\rangle \\{} \\
     \hskip5cm= \langle G_i(u^k_{1,j+\frac{1}{2}},\cdots,u^k_{4,j+\frac{1}{2}})\rangle,
	\end{array}
	\end{equation}
	which is consistent with the first equation of (\ref{Cross-Diffusion}).
\end{proposition}

\subsection{Boundary conditions }{\label{SUBBC}}
For the numerical solution of the kinetic equation (\ref{Cross-Diffusion}), the following inflow boundary conditions are usually prescribe for the distribution function $f$ 
$$f(t,x_{\min},v)=f_{l}(v),\quad v>0\qquad \hbox{and}\qquad f(t,x_{\max},v)=f_{r}(v),\quad v<0.$$
Moreover, the inflow boundary conditions can be rewritten in the micro-macro formulation \eqref{mM1} by
$$u_i(t,x_0)M+\frac{\varepsilon}{2}(g_i(t,x_{\frac{1}{2}},v)+g_i(t,x_{-\frac{1}{2}},v))=f_{i,l}(v),\quad v<0,$$
$$u_i(t,x_{N_x})M+\frac{\varepsilon}{2}(g_i(t,x_{N_x+\frac{1}{2}},v)+g_i(t,x_{N_x-\frac{1}{2}},v))=f_{i,r}(v),\quad v>0.$$
We consider the following artificial Neumann boundary conditions for the other velocities:
$$g_i(t,x_{\frac{1}{2}},v_l)=g_i(t,x_{-\frac{1}{2}},v_l),\quad v<0,$$
$$g_i(t,x_{N_x+\frac{1}{2}},v_l)=g_i(t,x_{N_x-\frac{1}{2}},v_l),\quad v>0.$$
Furthermore, the ghost points can be computed as follows:
\begin{equation}
g_{i,j-\frac{1}{2}}^{k+1}=\left\{
\begin{array}{l}
\frac{2}{\varepsilon}(f_{l}(v_\ell)-u_{i,0}^{k+1}M)-g_{i,\frac{1}{2},l}^{k+1},\quad v>0,  \\
{}  \\
g_{i,\frac{1}{2}}^{k+1},\quad v<0;
\end{array} \right.
\end{equation}
\begin{equation}
g_{i,N_x+\frac{1}{2}}^{k+1}=\left\{
\begin{array}{l}
\frac{2}{\varepsilon}(f_{r}(v_l)-u_{i,N_x}^{k+1}M)-g_{i,N_x-\frac{1}{2}}^{k+1},\quad v<0,  \\
{}  \\
g_{i,N_x-\frac{1}{2}}^{k+1},\quad v>0.
\end{array} \right. 
\end{equation}
Finally, from (\ref{mMDD2}) we have 
\begin{equation}
\left\{
\begin{array}{l}
\big(1+\frac{2\Delta t}{\varepsilon \Delta x}\langle v^+M_{i}\rangle\big)u_{i,0}^{k+1}=u_{i,0}^k-\frac{\Delta t}{\Delta x} \Big\langle(v+v^+-v^-)g_{i,\frac{1}{2}}^{k+1}-\frac{2v_l^+}{\varepsilon}f_{l}(v)\Big\rangle+\Delta t G_i(u_{1,0}^k,\cdots,u_{4,0}^k),\\
\\
{}  \\
\big(1-\frac{2\Delta t}{\varepsilon \Delta x}\langle v^-M_{i}\rangle\big)u_{i,N_x}^{k+1}=u_{i,N_x}^k-\frac{\Delta t}{\Delta x} \big\langle\frac{2v^-}{\varepsilon}f_{r}(v)-(v-v^++v^-)g_{i,N_x-\frac{1}{2}}^{k+1}\big\rangle+\Delta tG_i(u_{1,N_x}^k,\cdots,u_{4,N_x}^k). 
\end{array} \right. 
\end{equation}

\section{Numerical simulations} \label{Sec4}
In this section we provide some numerical simulations obtained from micro-macro formulation scheme presented in Section \ref{Sec3} and from the finite difference scheme of system \eqref{Cross-Diffusion}. Firstly, we show the asymptotic preserving scheme property, such as the uniform stability with respect to the parameter $\varepsilon$ and the consistence with the diffusion limit. Secondly, we demonstrate the effect of presence of the diffusion terms in system \eqref{Cross-Diffusion}, and its sensitivity with respect to the different choices of the reproduction ratio $R_0$. Finally, we attempt to model the COVID-19 pandemic in Morocco, where the data is available at \cite{[Rf2]} where we demonstrate the importance in considering a time-dependent rate transmission.\\
For the numerical simulations, the velocity space is $V = [-1,1]$ with number of grids $N_v = 164$, which can provide sufficient accuracy for numerical simulations (see, e.g., \cite{[KY13]}) and the time step is $t = 10^{-3}$, and we consider  the space domain $x = [-2, 2]$ with number of cells $N_x = 200$ and the periodic boundary condition. Moreover, we adopt a set of parameters as an example to analyse the results by varying some of them: $\mu=1/83,\; \alpha=0.06,\; \xi=1/4,\;\gamma=1/8$ as in \cite{[CSB20]}. For the diffusion coefficients, two cases are considered: i) case without diffusion and case of diffusion coefficients $d_1 = 0.05,\; d_2 = 0.025, \;d_3 = 0.001$, and $d_4 = 0$, similarly to \cite{[SSL10]}. Finally, we use the following two different cases of initial conditions:
\begin{itemize}
	\item[$i)$]
	\begin{equation}
	\left\{
	\begin{array}{l}
	S_0=2.6\,\Big(\exp(-(\frac{x-0.5}{0.12})^2)+\exp(-(\frac{x+0.5}{0.12})^2)\Big)/(0.9\,\pi),\\
   I_0=0.04\,\exp(-2\,x^2),	\\
	E_0=R_0=0, \\
	N_0=S_0+I_0, 
	\end{array} \right. 
	\end{equation}
	
	\item[$ii)$] 	
	\begin{equation}
	\left\{
	\begin{array}{l}
	S_0=0.96\,\exp(-10(\frac{x}{1.4})^2),\\
	I_0=0.04\,\exp(-2\,x^2),	\\
	E_0=R_0=0, \\
	N_0=S_0+I_0. 
	\end{array} \right. 
	\end{equation}
\end{itemize}
%\subsection{AP scheme}

\begin{figure}
	%\begin{minipage}{0.99\linewidth}
		\centering
	\begin{tikzpicture}[spy using outlines={ ,blue,magnification=3,size=1.5cm, connect spies}]
		\node{\includegraphics[height=1.7in ,width=1.8in]{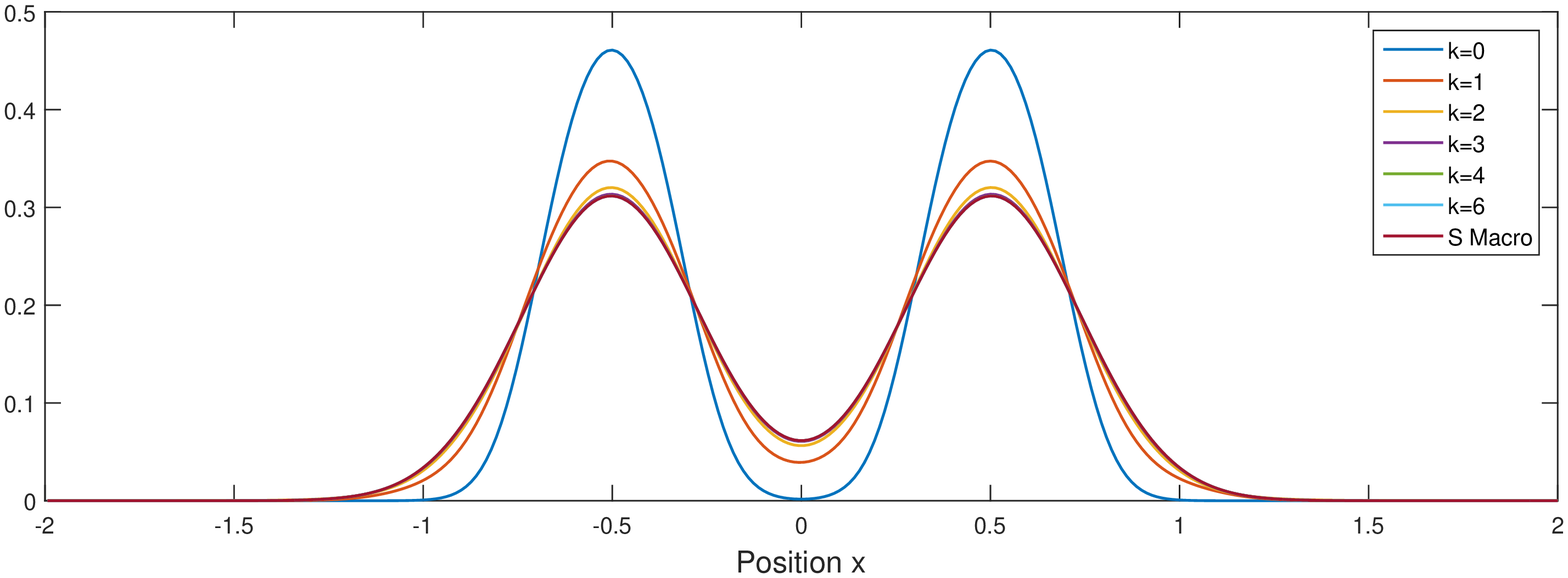}};
		\spy on (0.1,-1.2) in node [right] at  (.5,-1.45);
		\end{tikzpicture}
	\begin{tikzpicture}[spy using outlines={ ,blue,magnification=3,size=1.5cm, connect spies}]
	\node{\includegraphics[height=1.7in ,width=1.8in]{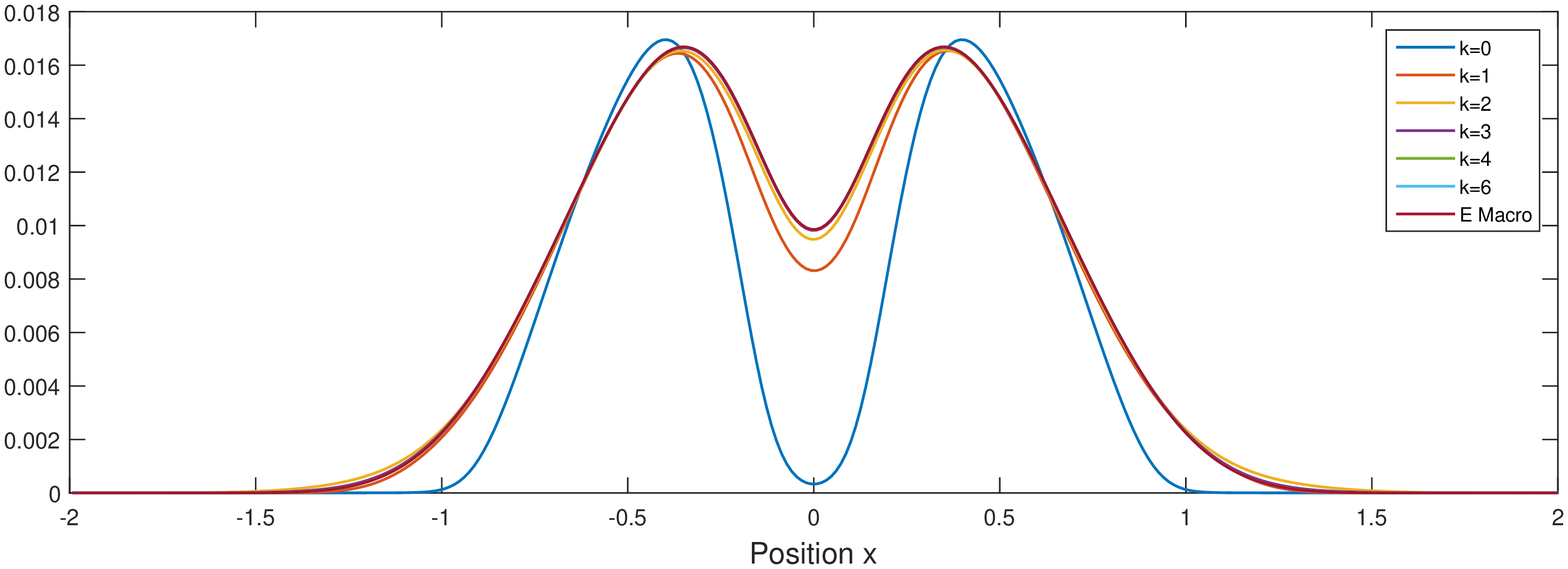}};
	\spy on (0.1,0.4) in node [right] at  (.5,-1.45);
	\end{tikzpicture}
	\begin{tikzpicture}[spy using outlines={ ,blue,magnification=3,size=1.5cm, connect spies}]
	\node{\includegraphics[height=1.7in ,width=1.8in]{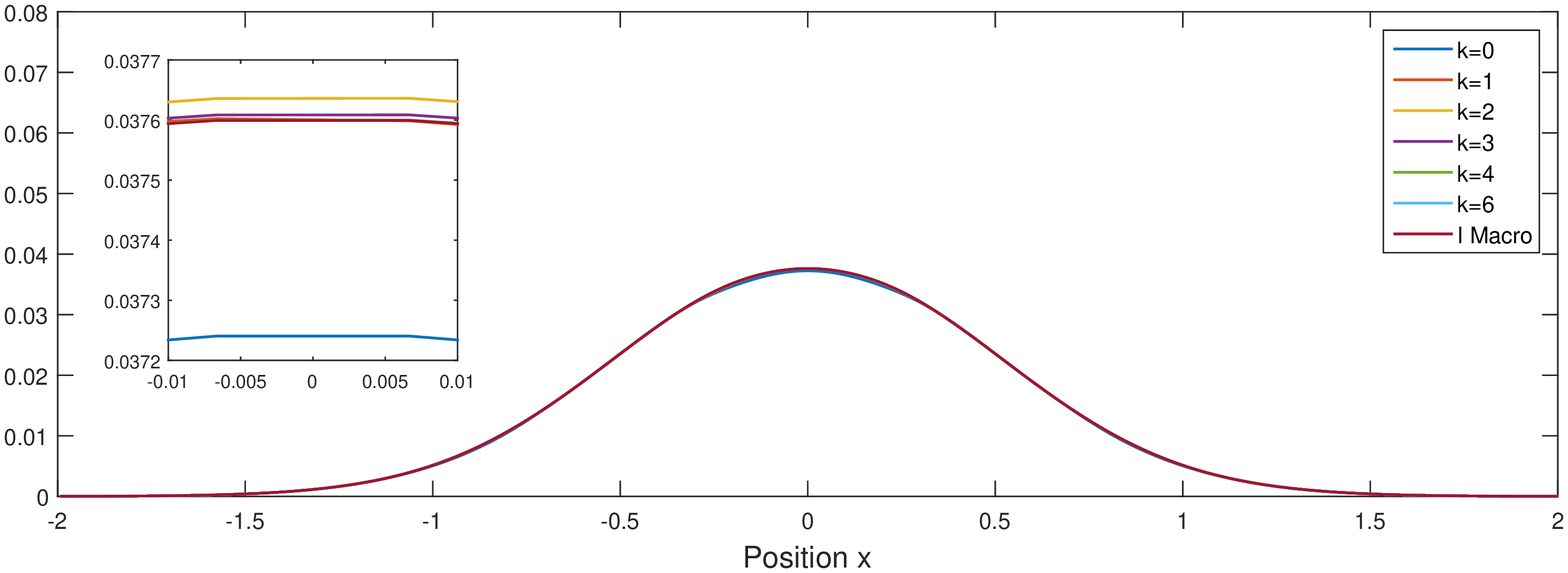}};
	\spy on (0.1,0.) in node [right] at  (.5,-1.45);
	\end{tikzpicture}
	
	\begin{tikzpicture}[spy using outlines={ ,blue,magnification=3,size=1.5cm, connect spies}]
	\node{\includegraphics[height=1.7in ,width=1.8in]{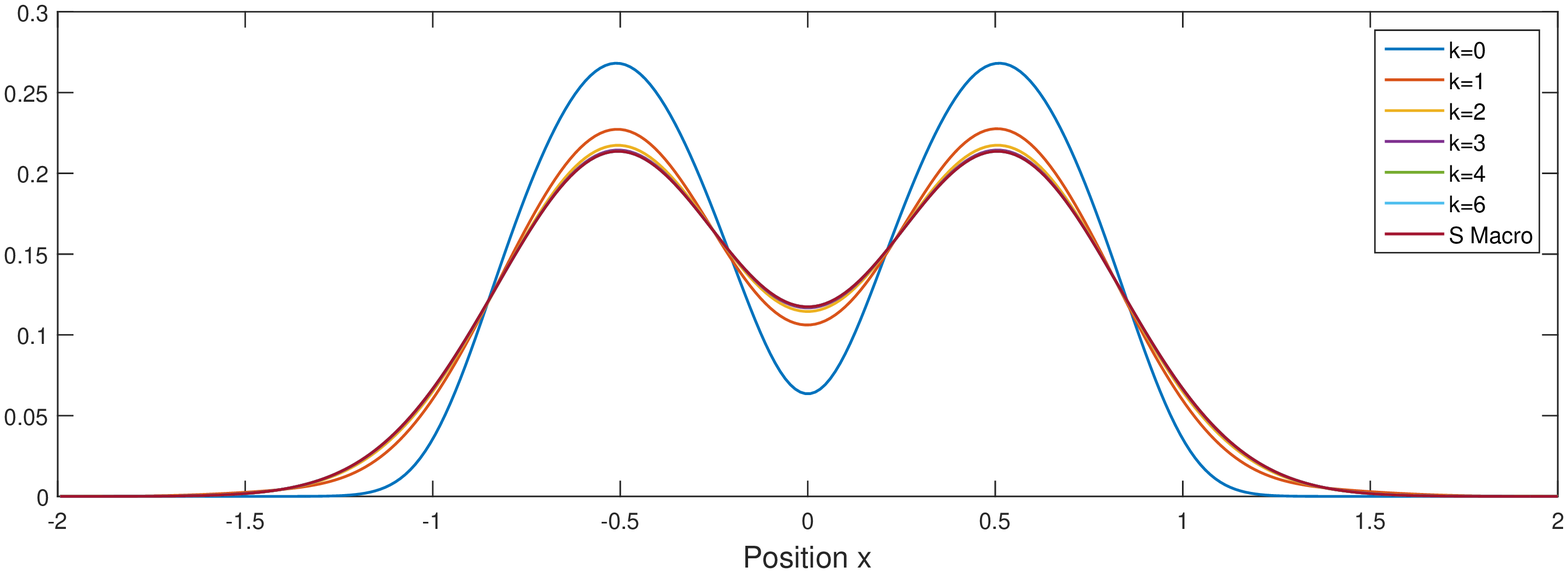}};
	\spy on (0.1,-0.2) in node [right] at  (.5,-1.45);
	\end{tikzpicture}
	\begin{tikzpicture}[spy using outlines={ ,blue,magnification=3,size=1.5cm, connect spies}]
	\node{\includegraphics[height=1.7in ,width=1.8in]{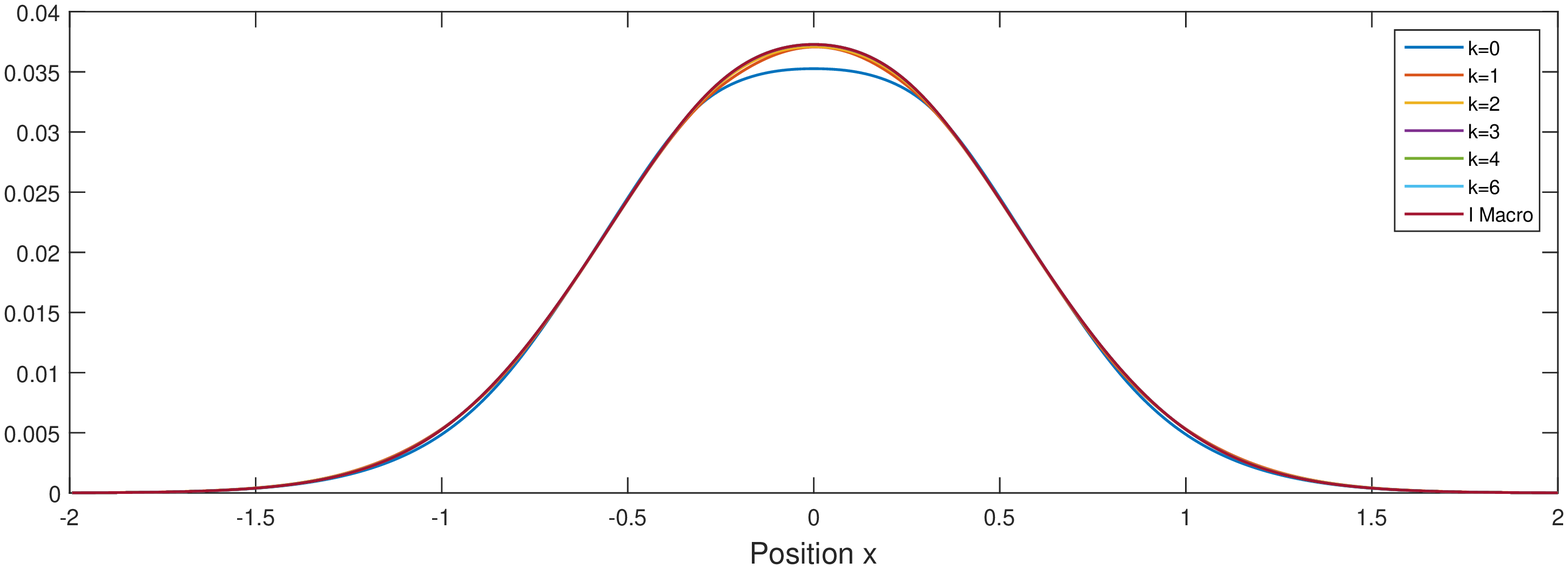}};
	\spy on (0.1,1.4) in node [right] at  (.5,-1.45);
	\end{tikzpicture}
	\begin{tikzpicture}[spy using outlines={ ,blue,magnification=3,size=1.5cm, connect spies}]
	\node{\includegraphics[height=1.7in ,width=1.8in]{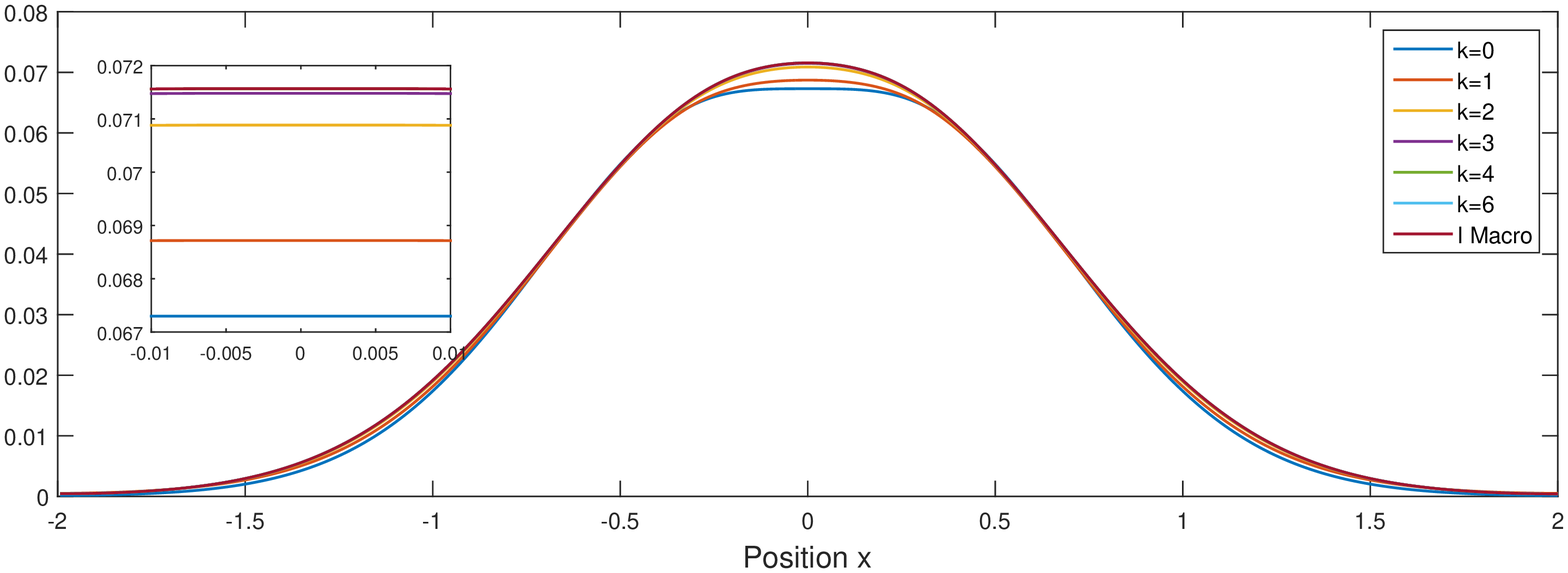}};
	\spy on (0.1,1.4) in node [right] at  (.5,-1.45);
	\end{tikzpicture}
	
		\begin{tikzpicture}[spy using outlines={ ,blue,magnification=3,size=1.5cm, connect spies}]
	\node{\includegraphics[height=1.7in ,width=1.8in]{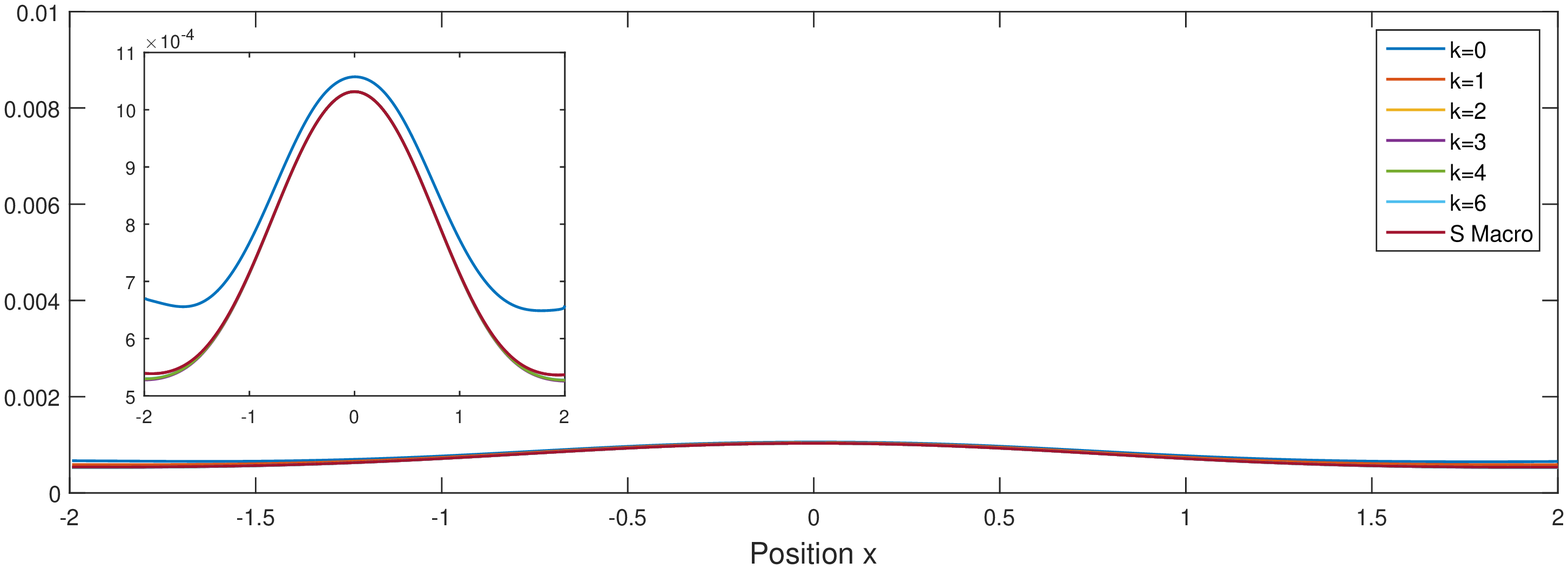}};
	\spy on (0.1,-1.2) in node [right] at  (.5,-1.45);
	\end{tikzpicture}
	\begin{tikzpicture}[spy using outlines={ ,blue,magnification=3,size=1.5cm, connect spies}]
	\node{\includegraphics[height=1.7in ,width=1.8in]{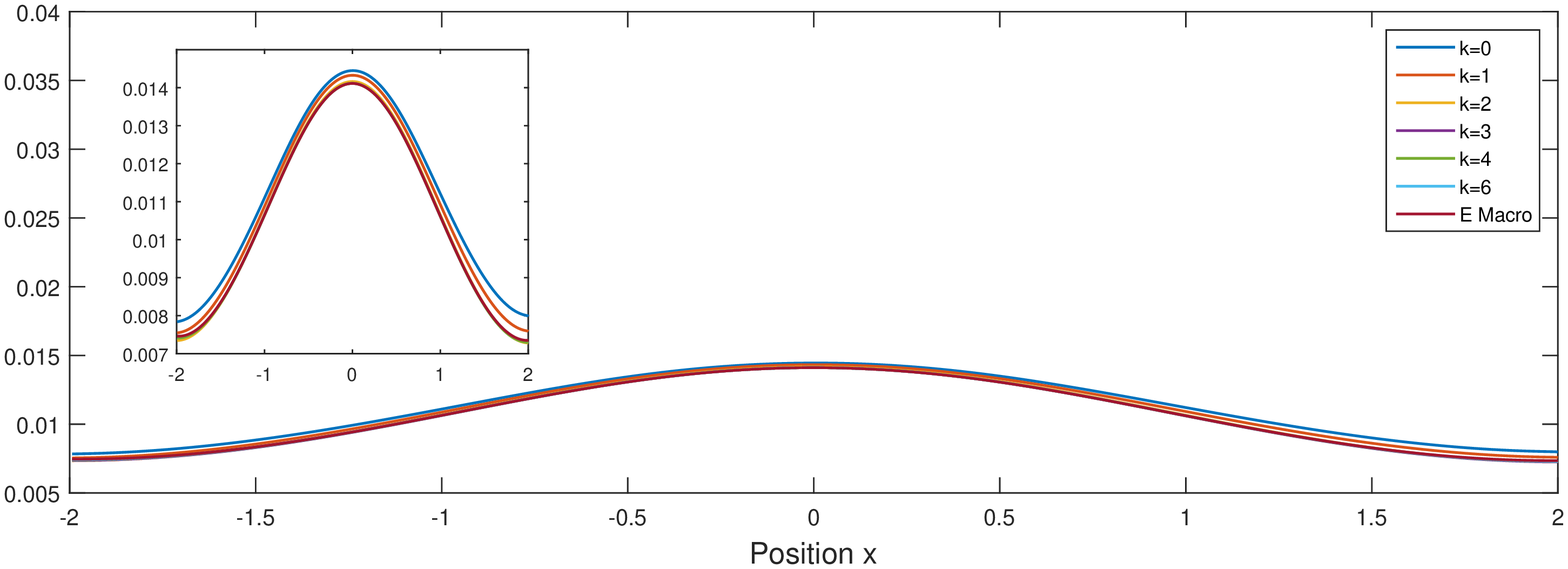}};
	\spy on (0.1,-0.7) in node [right] at  (.5,-1.45);
	\end{tikzpicture}
	\begin{tikzpicture}[spy using outlines={ ,blue,magnification=3,size=1.5cm, connect spies}]
	\node{\includegraphics[height=1.7in ,width=1.8in]{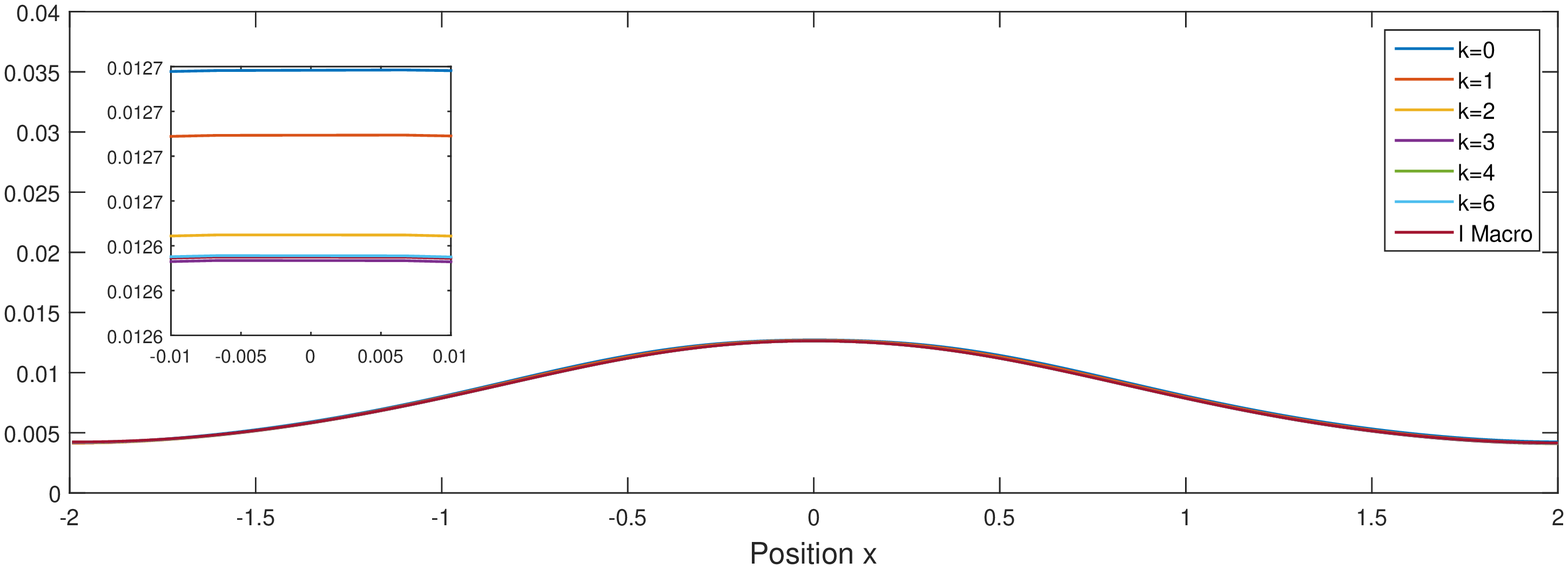}};
	\spy on (0.1,-.5) in node [right] at  (.5,-1.45);
	\end{tikzpicture}
	
		\begin{tikzpicture}[spy using outlines={ ,blue,magnification=3,size=1.5cm, connect spies}]
	\node{\includegraphics[height=1.7in ,width=1.8in]{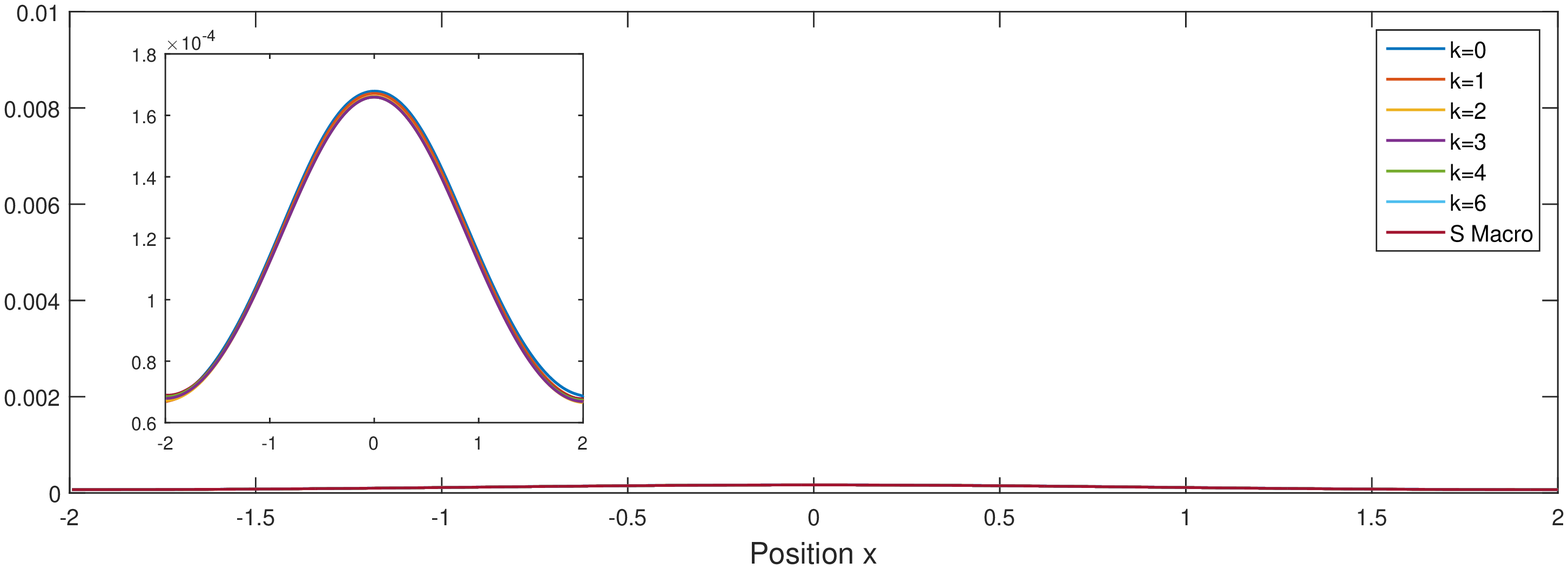}};
	\spy on (0.1,-1.6) in node [right] at  (.5,-1.45);
	\end{tikzpicture}
	\begin{tikzpicture}[spy using outlines={ ,blue,magnification=3,size=1.5cm, connect spies}]
	\node{\includegraphics[height=1.7in ,width=1.8in]{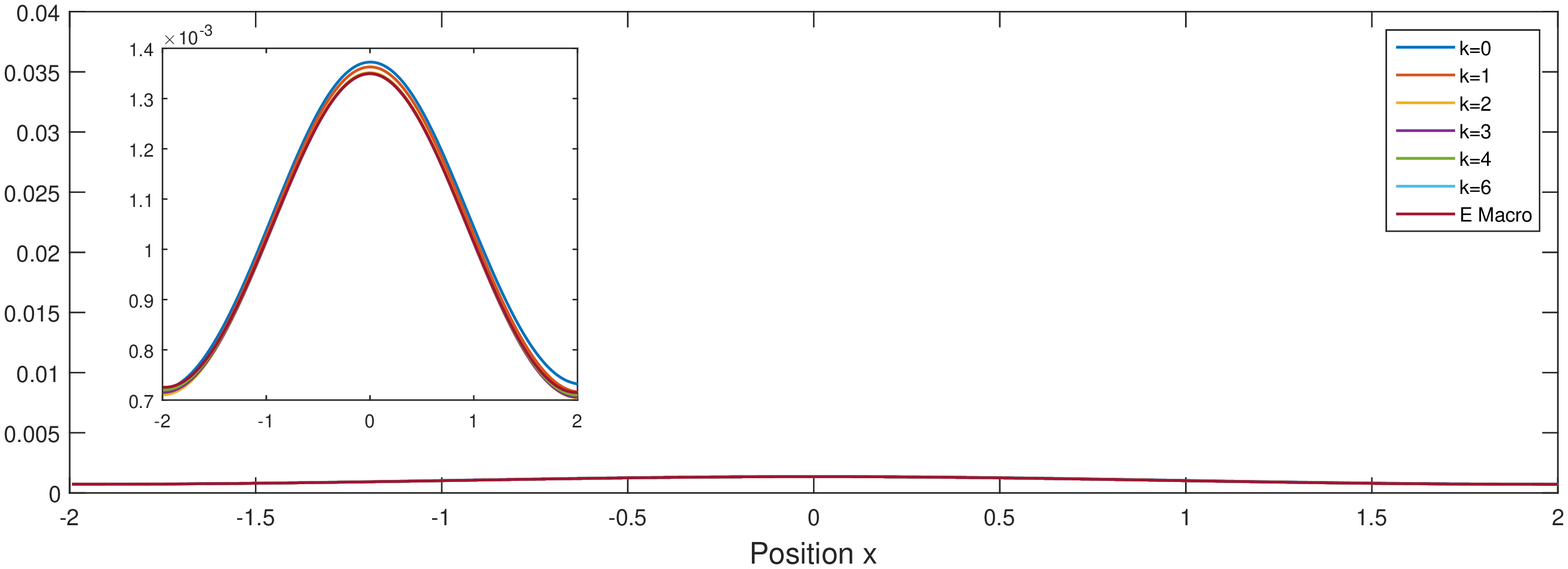}};
	\spy on (0.1,-1.6) in node [right] at  (.5,-1.45);
	\end{tikzpicture}
	\begin{tikzpicture}[spy using outlines={ ,blue,magnification=3,size=1.5cm, connect spies}]
	\node{\includegraphics[height=1.7in ,width=1.8in]{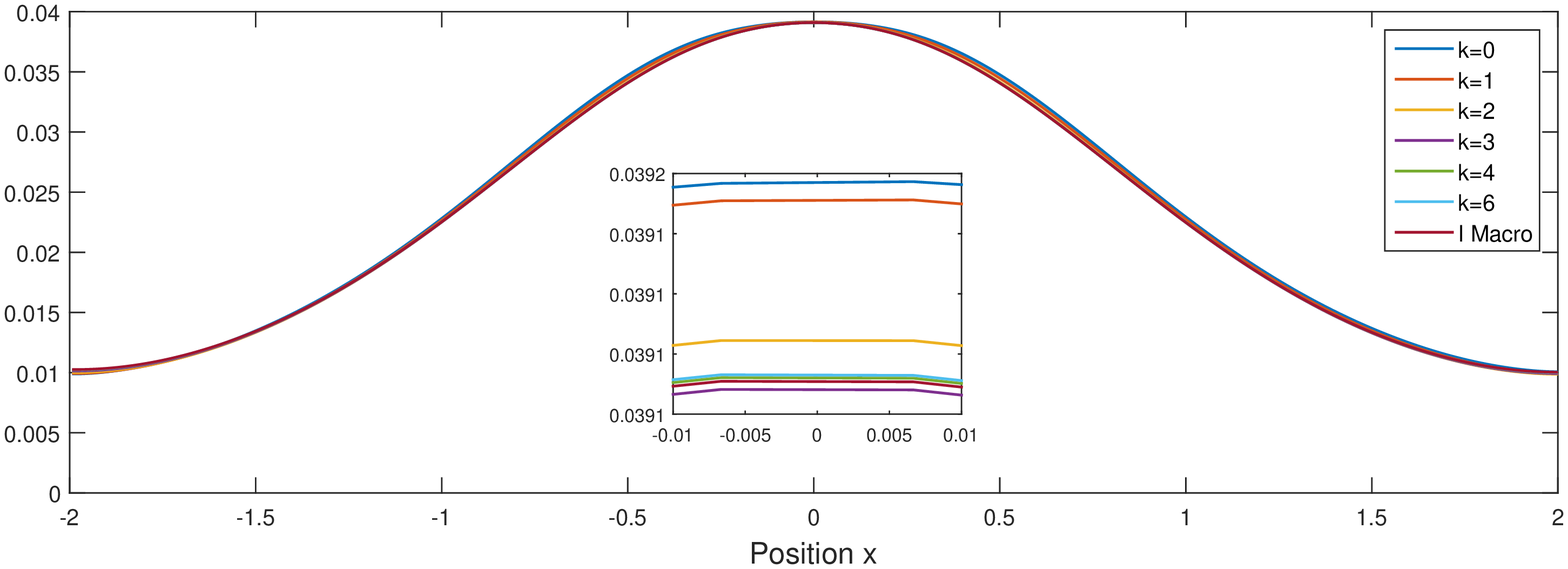}};
	\spy on (0.1,1.7) in node [right] at  (.5,-1.45);
	\end{tikzpicture}	
		%\hfil
	%\end{minipage}
	\caption{Dynamics of the densities $S$ (first column), $E$ (second column) and $I$ (third column) obtained from the (AP)-scheme with $\varepsilon =2\times10^{-k}$, $k = 0,\,1,\, 2,\,3,\, 4,\, 6$ against the macroscopic scheme with initial conditions $i)$ at successive time $t=0.5,\,1,\,5,\,10$.}\label{F0}	
\end{figure}

\begin{figure}[h!]
	\centering
	\subfigure[]{\includegraphics[height=1.7in ,width=2in]{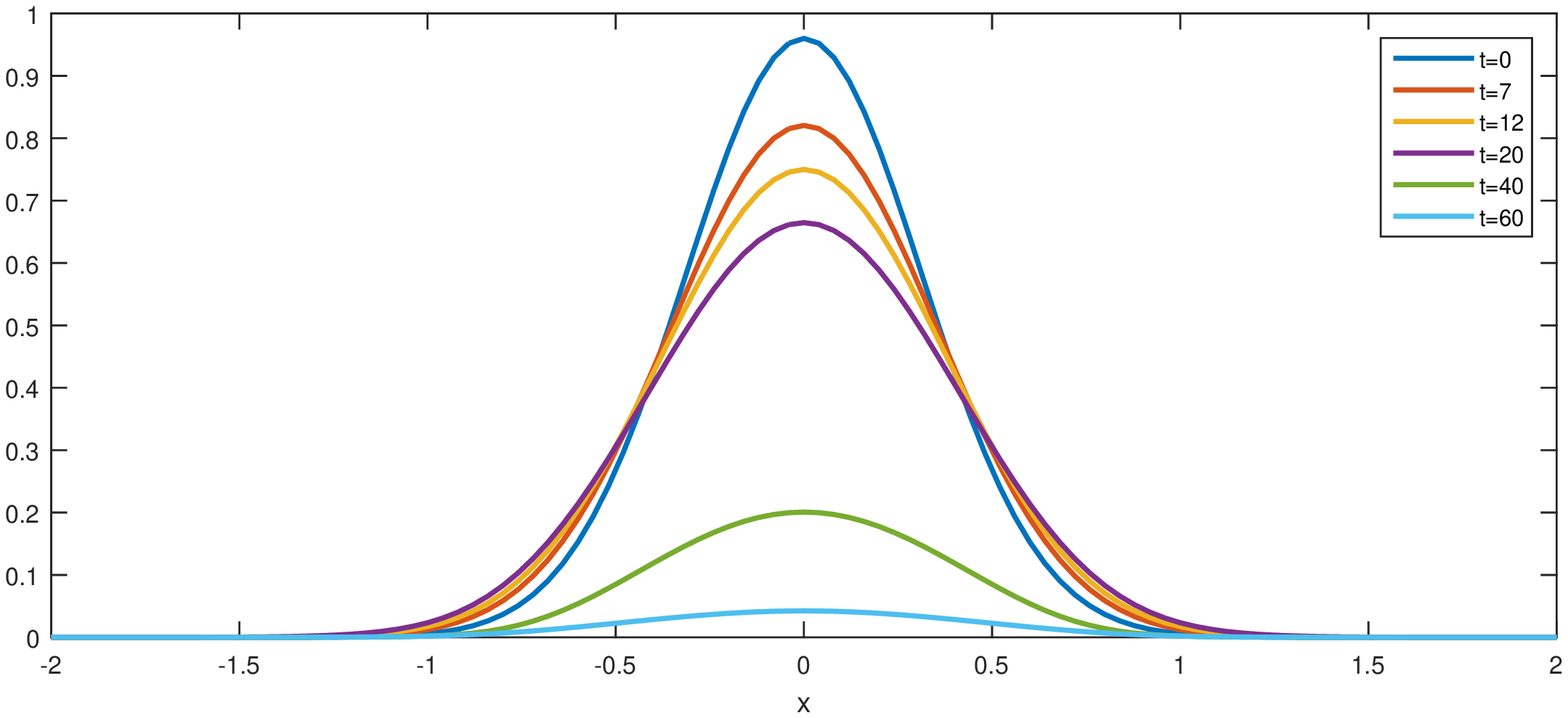}}~
	\subfigure[]{\includegraphics[height=1.7in ,width=2in]{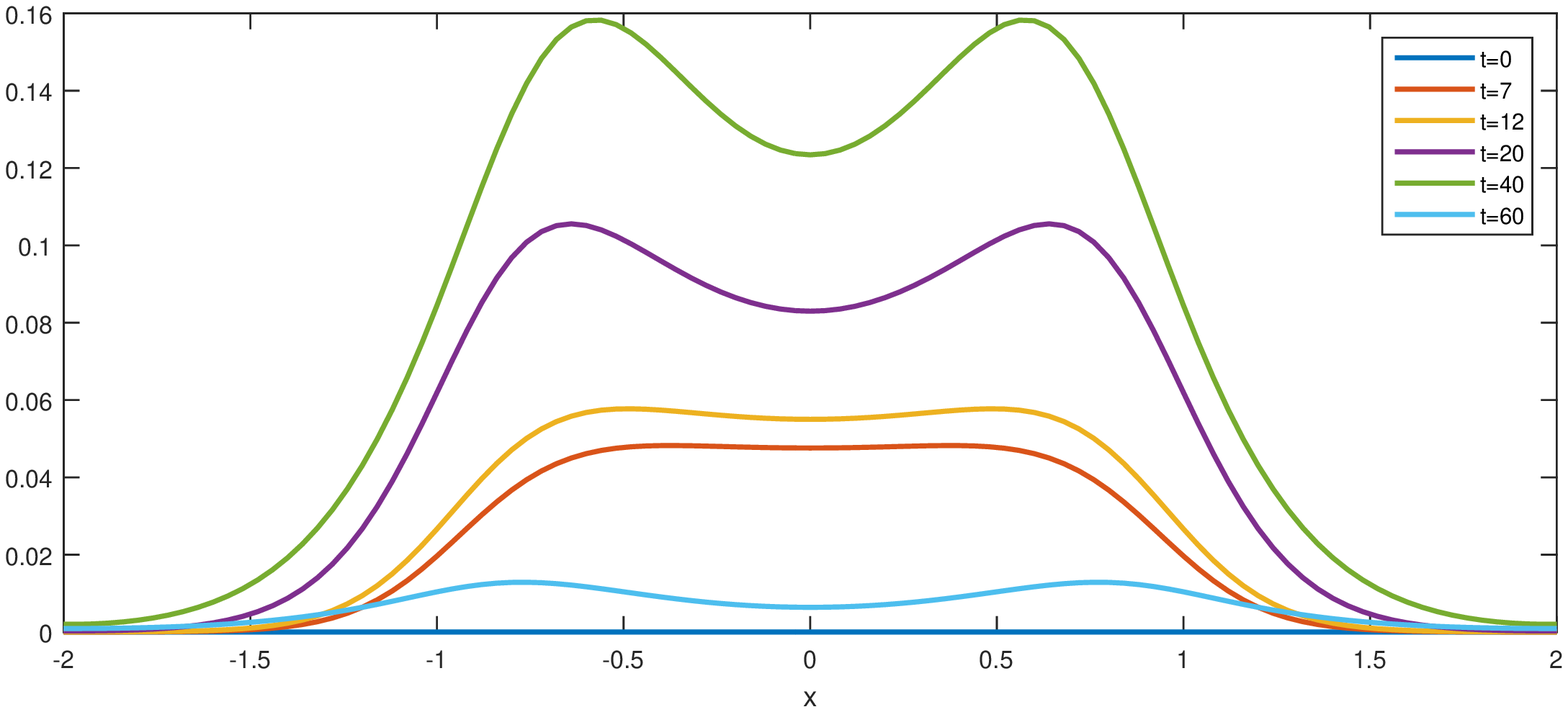}}~
	\subfigure[]{\includegraphics[height=1.7in ,width=2in]{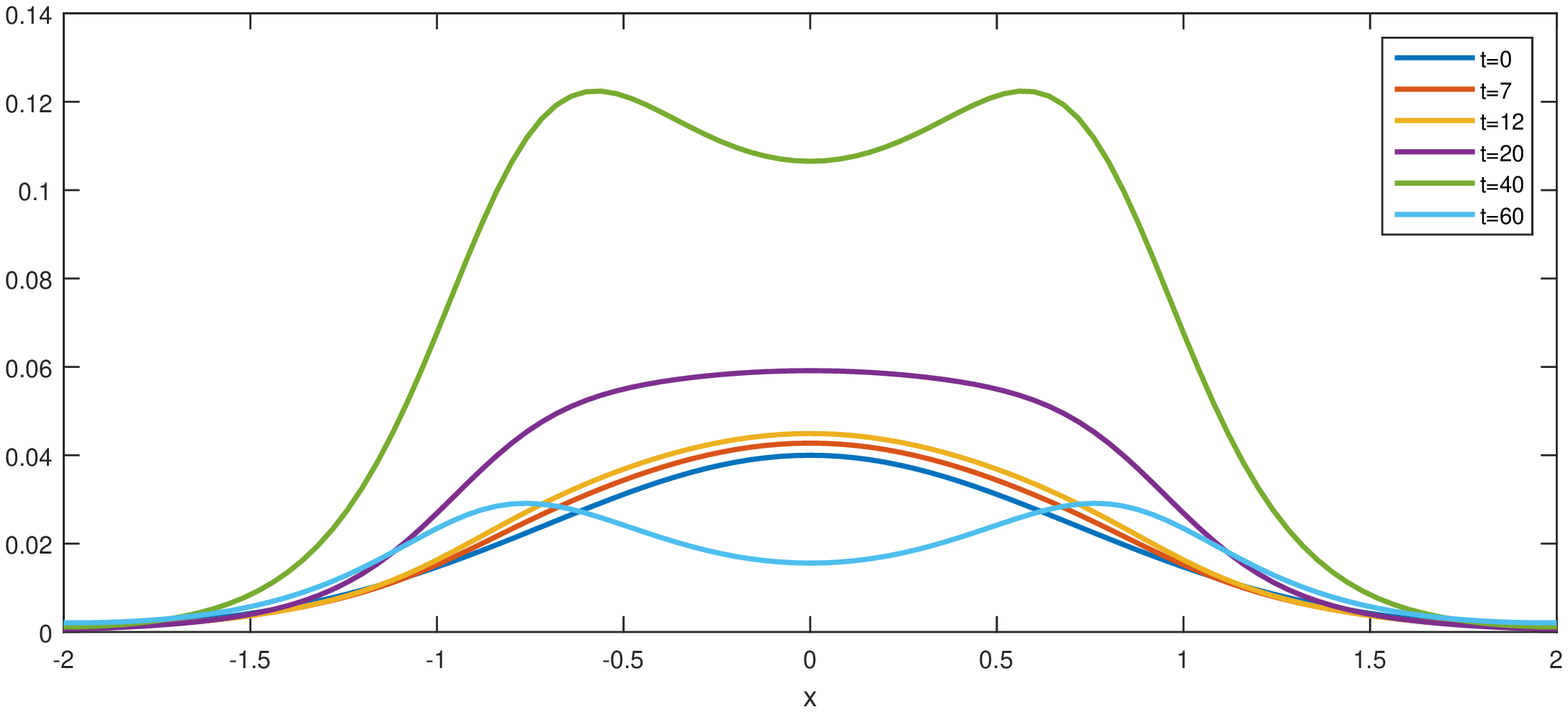}}
	\subfigure[]{\includegraphics[height=1.7in ,width=2in]{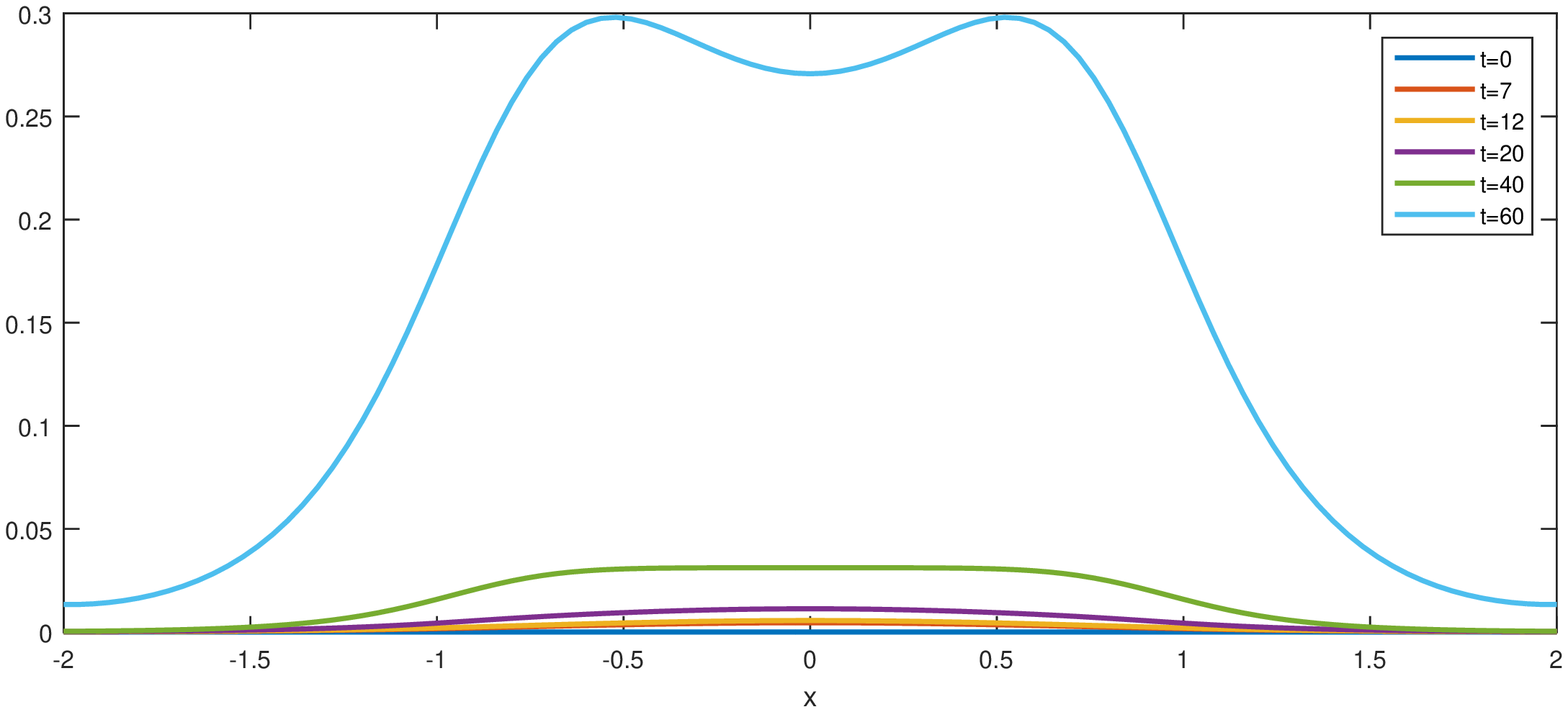}}~
	\subfigure[]{\includegraphics[height=1.7in ,width=2in]{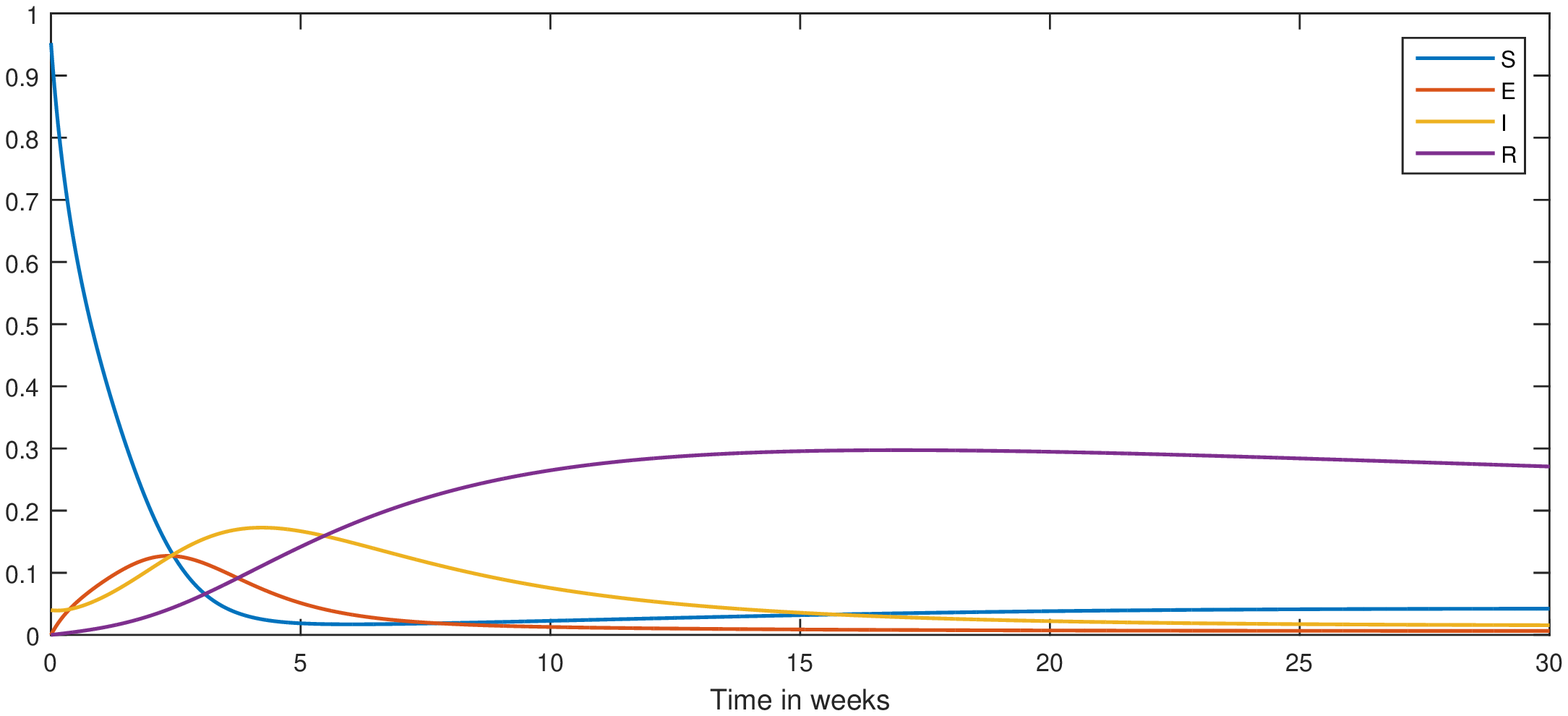}}~
	\subfigure[]{\includegraphics[height=1.7in ,width=2in]{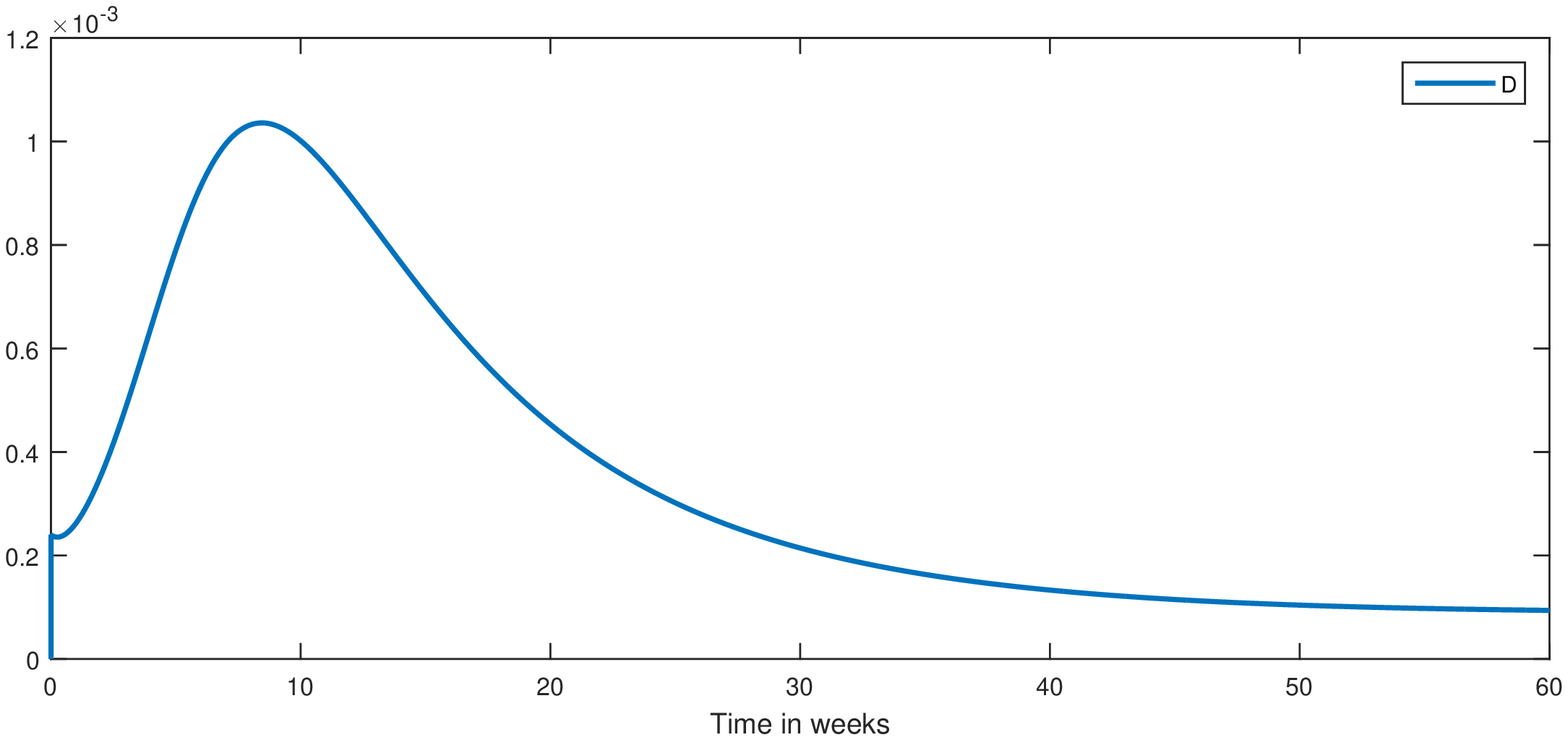}}
	\caption{Snapshot of the obtained numerical solutions of $S$ (sub-figure (a)), $E$ (sub-figure (b)), $I$ (sub-figure (c)), $R$ (sub-figure (d)) at successive time $t=0,\,7,\,12,\,30,\,40,\,60$ from (AP)-scheme with $\varepsilon=10^{-6}$ and its spatial variation at $x=0$ (sub-figure (e)), while (sub-figure (f)) is for the population of died individuals for the reproduction ratio $R_0=2$. }
	\label{F1}
\end{figure}

In Figure \ref{F0}, we present the numerical results obtained of susceptible, exposed and infected individuals from the (AP)-scheme and from the reaction-diffusion scheme with initial conditions $i)$ at successive times $t =  0.5,\, 1,\, 5,\,10$. We observe that the results obtained from the two schemes have almost the same profiles in the limit when the parameter $ \varepsilon= 2\times10^{-k}$, with $k=0,\,1,\,2,\,3,\,4,\,6$ goes to zero. This confirms that the asymptotic preserving scheme is uniformly stable along the transition from the kinetic regime to the macroscopic regime. Furthermore, we see that the asymptotic-preserving scheme converges better in time (see the zoomed  windows). Figure \ref{F1} shows the snapshot of the obtained numerical solutions from (AP)-scheme with $\varepsilon=10^{-6}$ of $S$ (sub-figure (a)), $E$ (sub-figure (b)), $I$ (sub-figure (c)), $R$ (sub-figure (d)) at successive time $t=0,\,7,\,12,\,30,\,40,\,60$ and their spatial variation at $x=0$ (sub-figure (e)), while (sub-figure (f)) is for the population of died individuals for the reproduction ratio $R_0=2$. \\

\begin{figure}[h!]
	\centering
	%\subfigure[]{\includegraphics[height=1.7in ,width=2.2in]{SN.eps}}~
	%\subfigure[]{\includegraphics[height=1.7in ,width=2.2in]{EN.eps}}~
	%\subfigure[]{\includegraphics[height=1.7in ,width=2.2in]{IN.eps}}
	\subfigure[]{\includegraphics[height=1.7in ,width=2.2in]{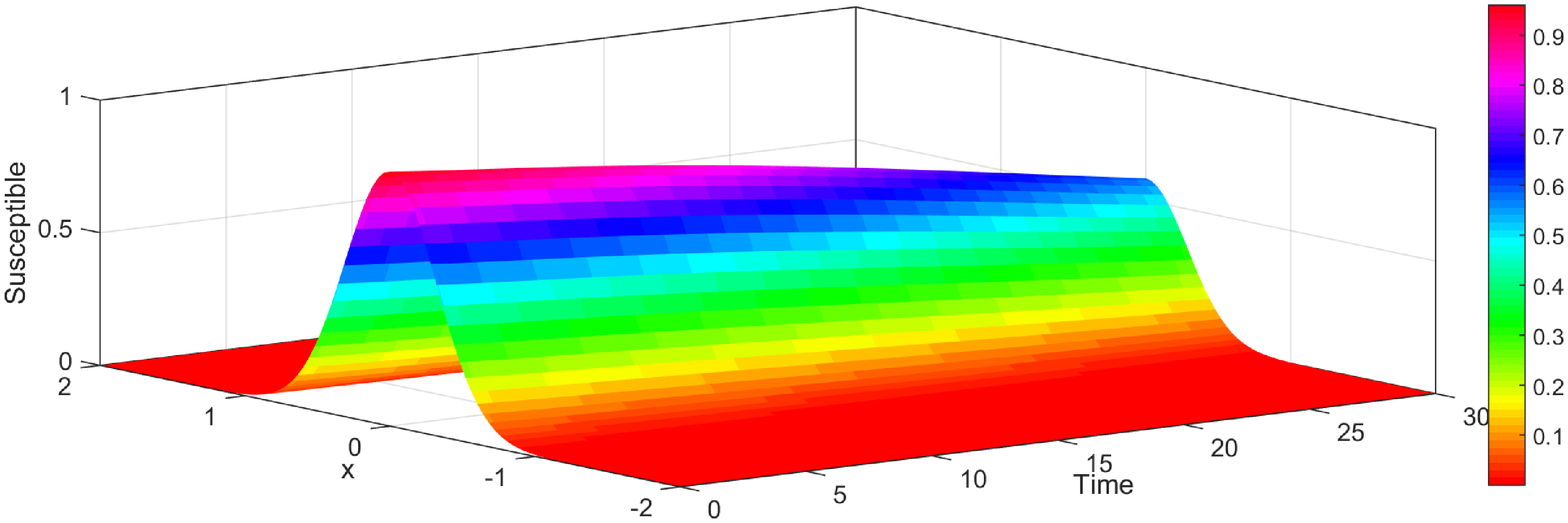}}~
	\subfigure[]{\includegraphics[height=1.7in ,width=2.2in]{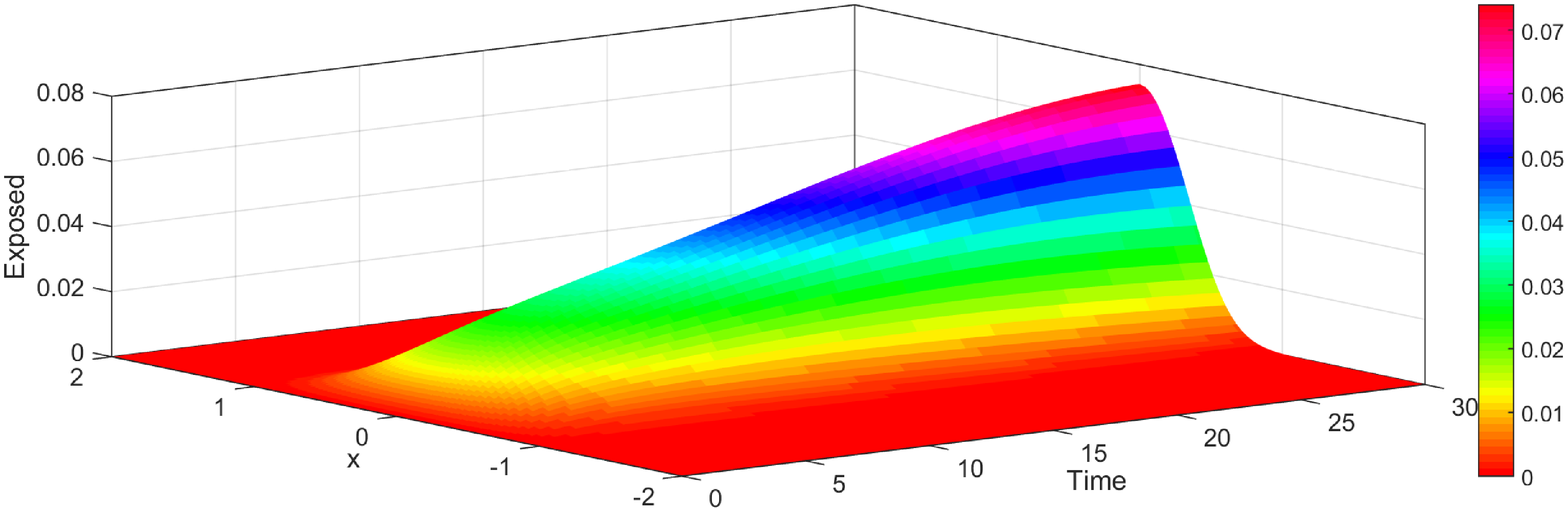}}~
	\subfigure[]{\includegraphics[height=1.7in ,width=2.2in]{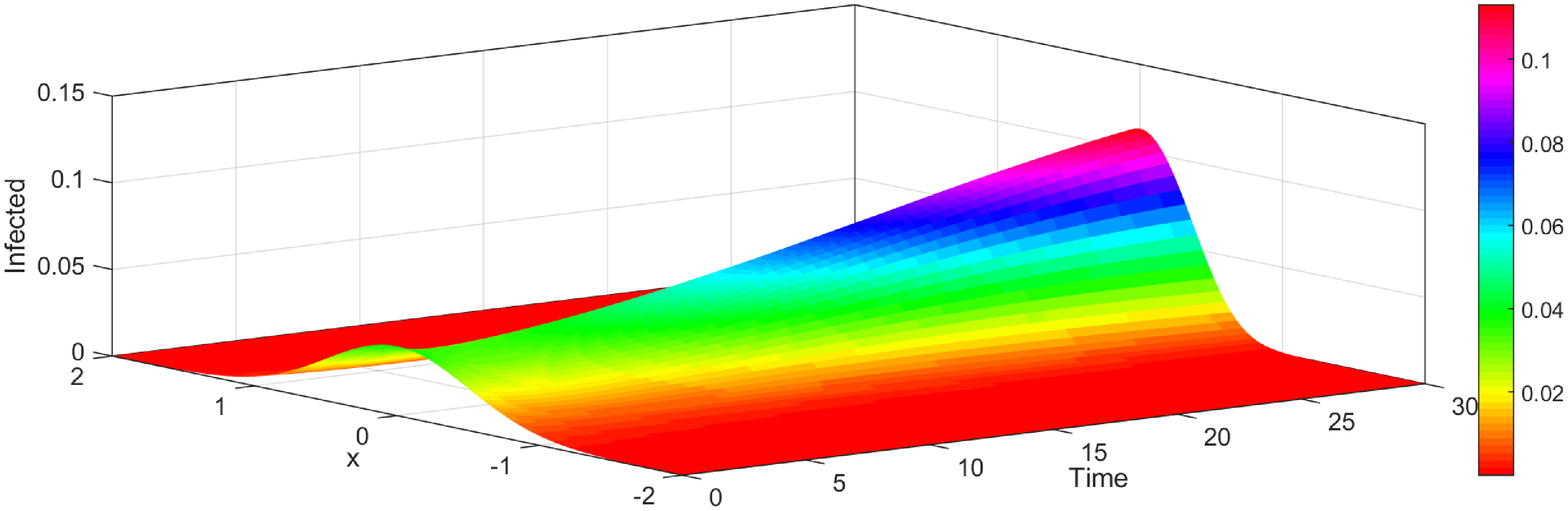}}
	%\subfigure[]{\includegraphics[height=1.7in ,width=2.2in]{SW.eps}}~
	%\subfigure[]{\includegraphics[height=1.7in ,width=2.2in]{EW.eps}}~
	%\subfigure[]{\includegraphics[height=1.7in ,width=2.2in]{IW.eps}}
	\subfigure[]{\includegraphics[height=1.7in ,width=2.2in]{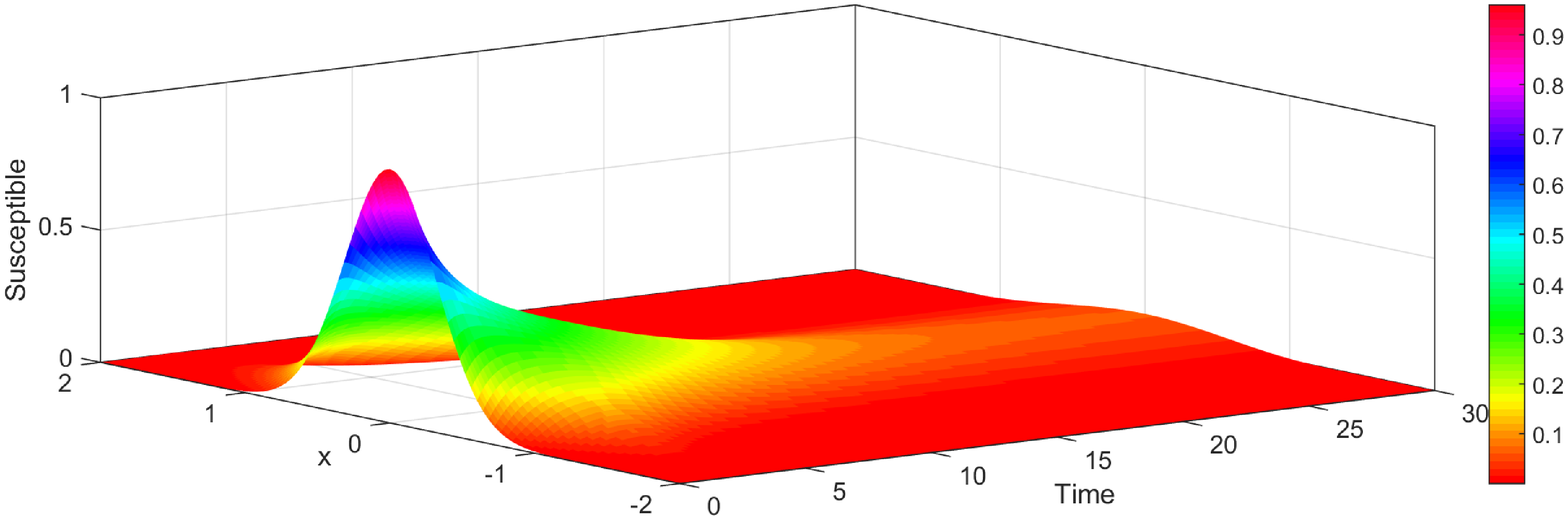}}~
	\subfigure[]{\includegraphics[height=1.7in ,width=2.2in]{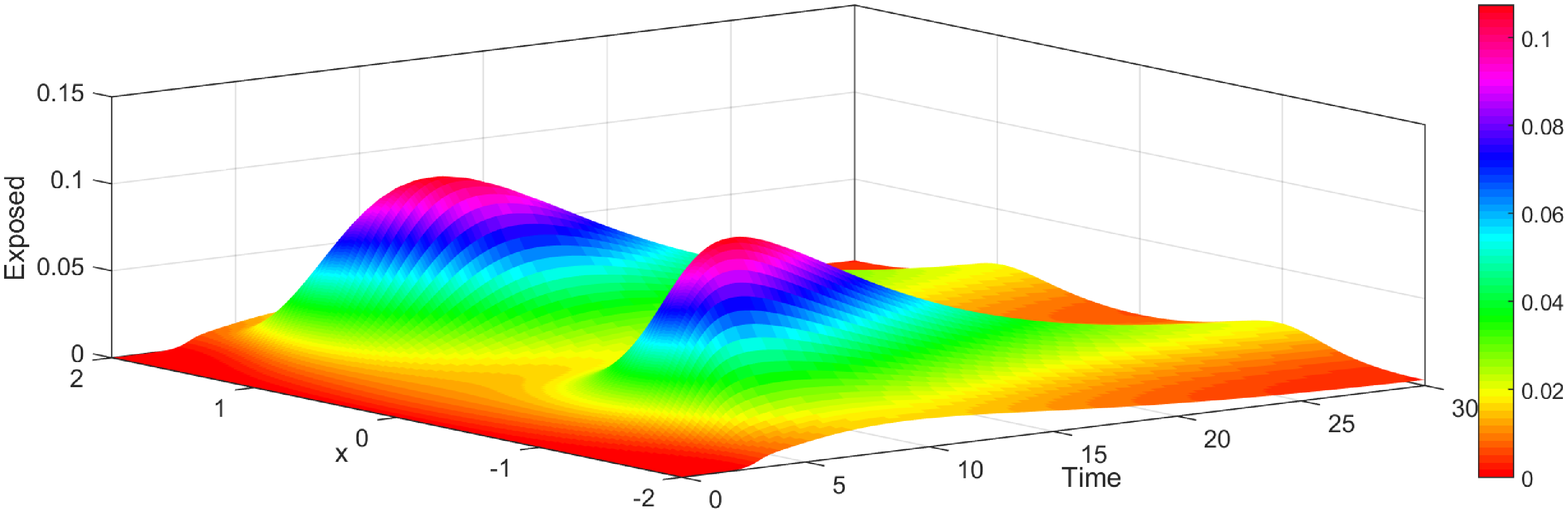}}~
	\subfigure[]{\includegraphics[height=1.7in ,width=2.2in]{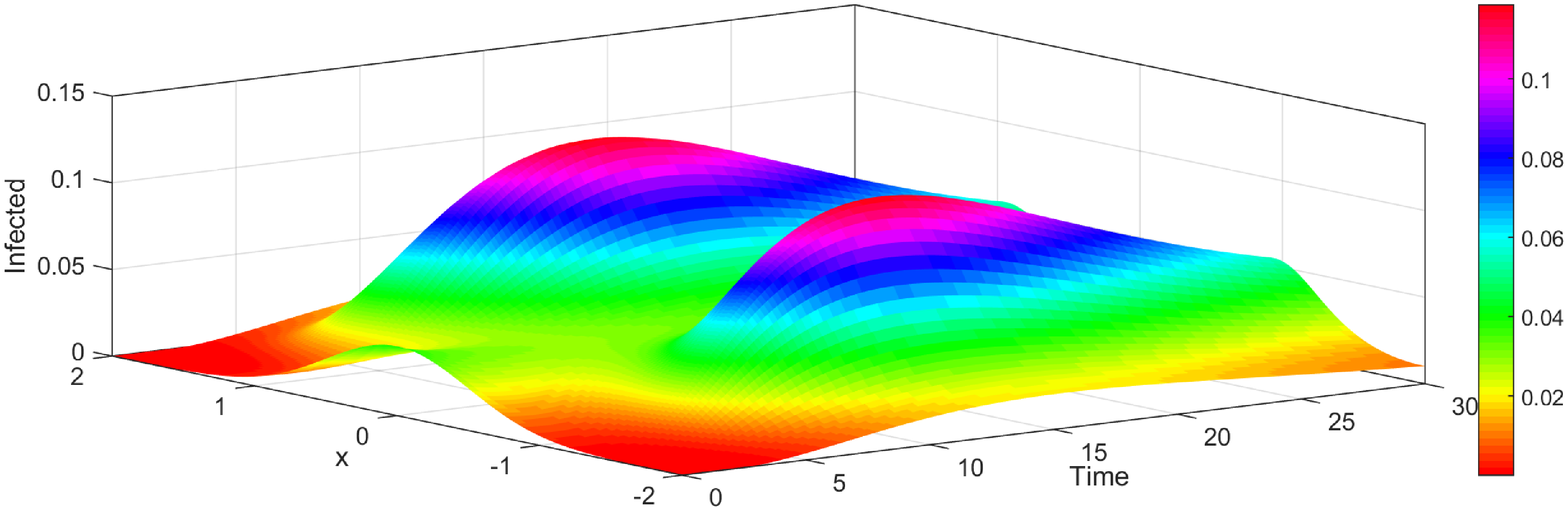}}
	\caption{Evolution of densities  $f_1,\,f_2,f_3$ using micro-macro scheme for $\varepsilon=10^{-6}$ and initial condition $ii)$ with diffusion (sub-figures (a), (b), (c)) and without diffusion (sub-figures (d), (e), (f)), respectively, for the reproduction ratio $R_0 = 2$.}
	\label{F3}
\end{figure}

To demonstrate the diffusion effect on individuals interactions, we consider the initial conditions $ii)$ with the reproduction ratio value $R_0=2$. Figure \ref{F3} provides the obtained results of susceptible, exposed and infected individuals from the (AP)-scheme with $\varepsilon=10^{-6}$ in the case without diffusion (see sub-figures (a)-(b)-(c)) and with diffusion (see sub-figures (d)-(e)-(f)).  We observe that the individuals are all centred around $x = 0$ in the absence of diffusion ($d_i=0$). While when the diffusion is considered, it is clear that the individuals are more spread along the $x-$axis.\\

\begin{figure}[h!]
	\centering
	\subfigure[]{\includegraphics[height=1.7in ,width=2in]{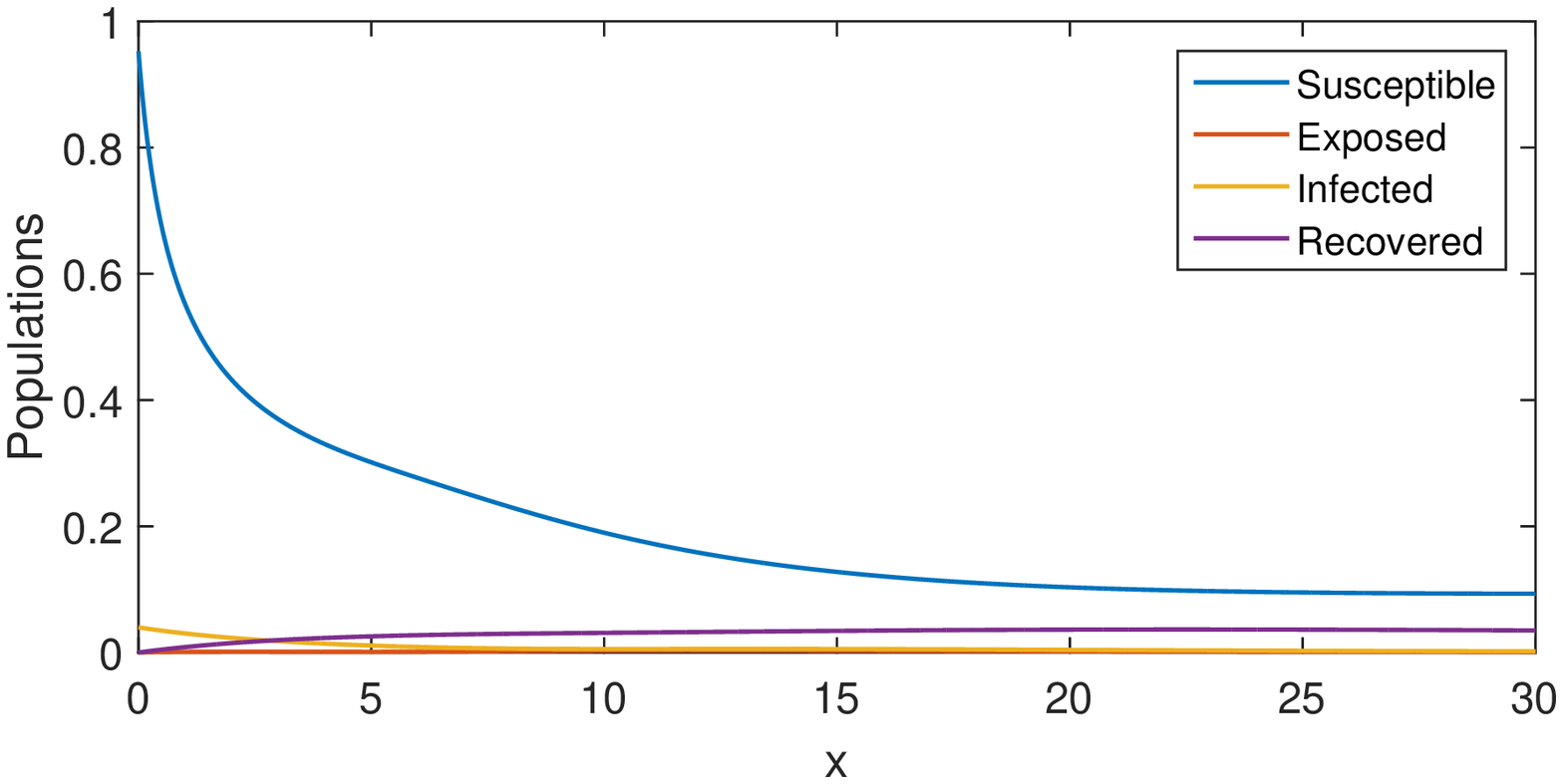}}~
	\subfigure[]{\includegraphics[height=1.7in ,width=2in]{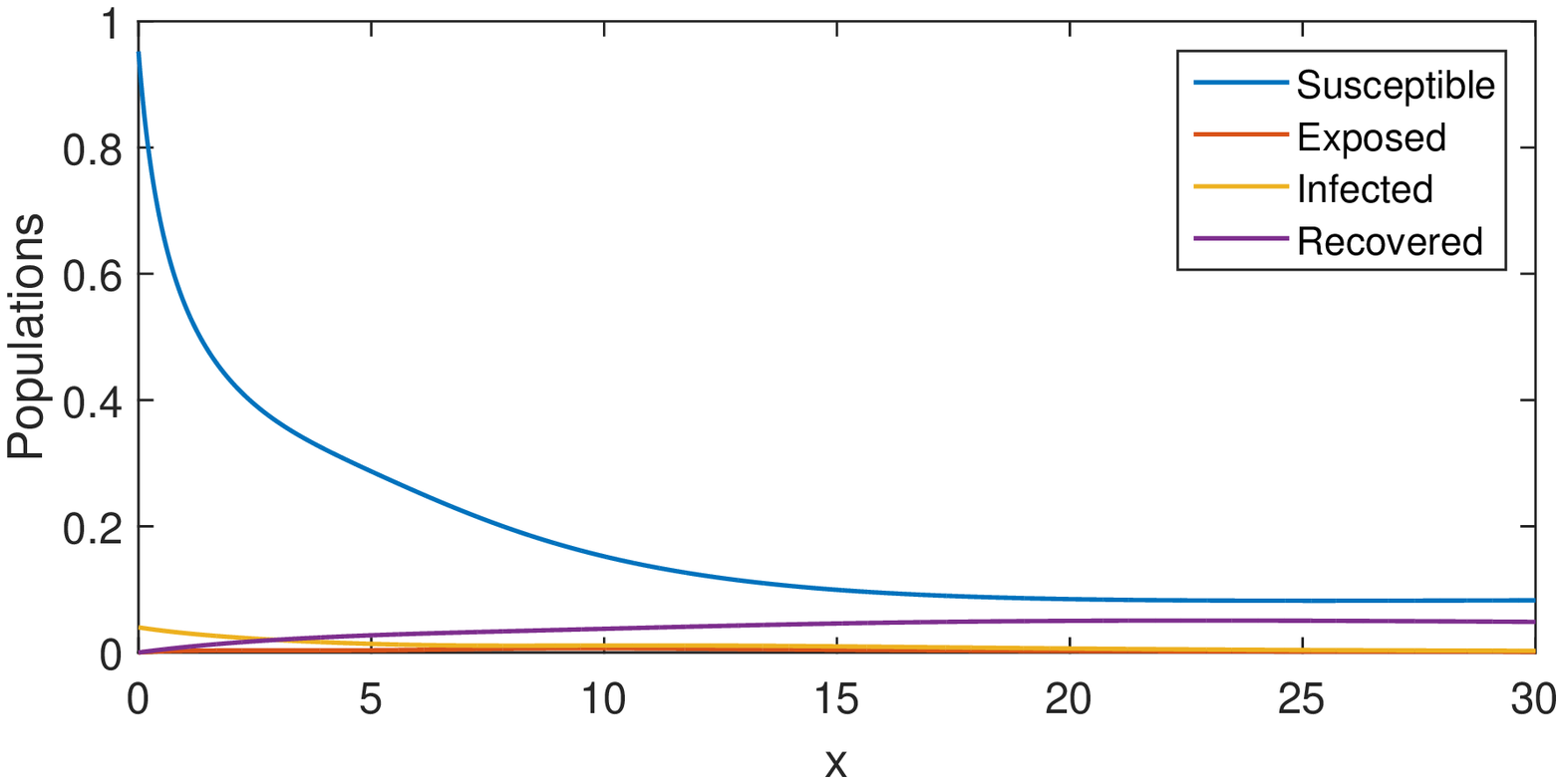}}~
	\subfigure[]{\includegraphics[height=1.7in ,width=2in]{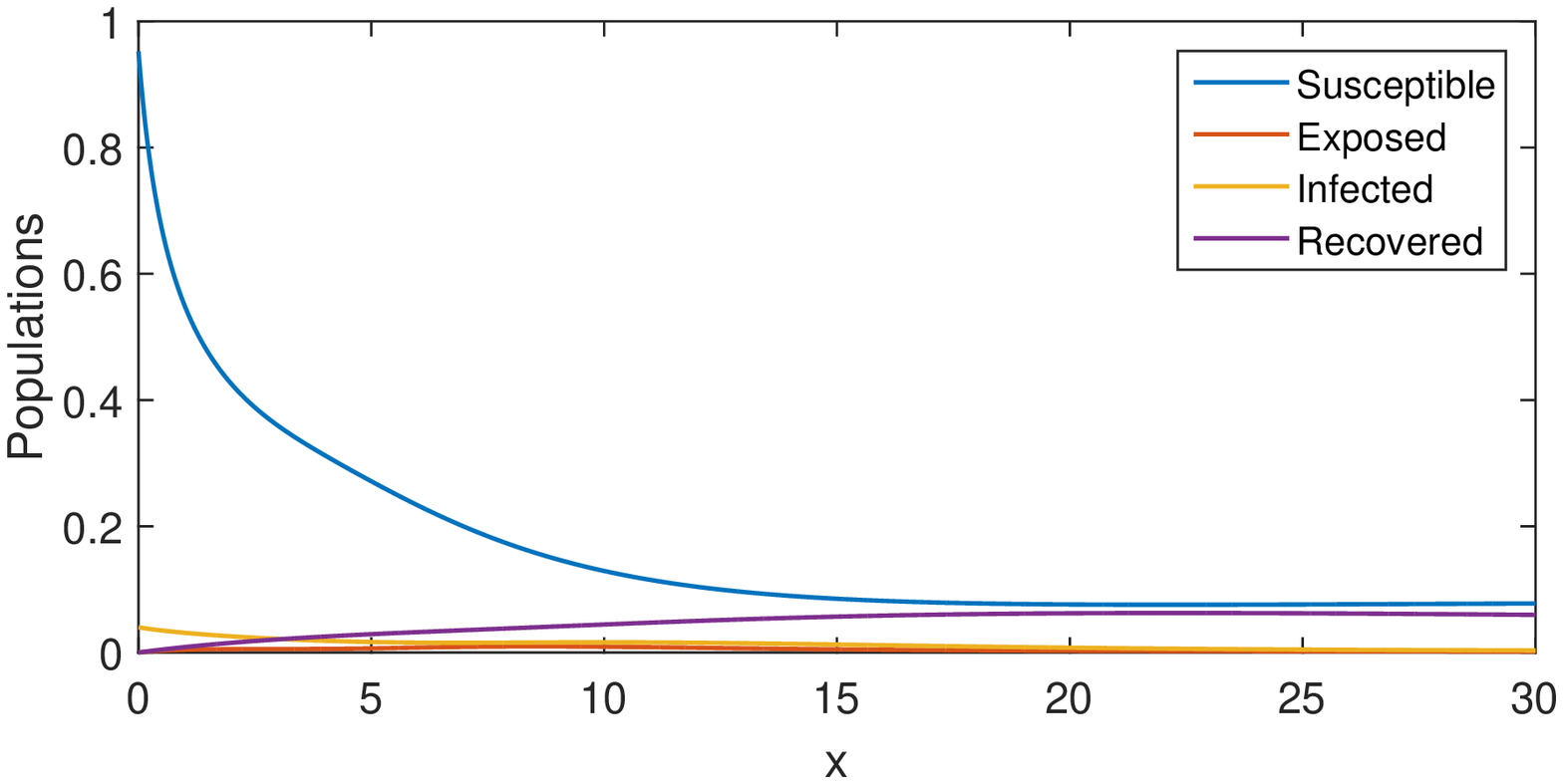}}
	\subfigure[]{\includegraphics[height=1.7in ,width=2in]{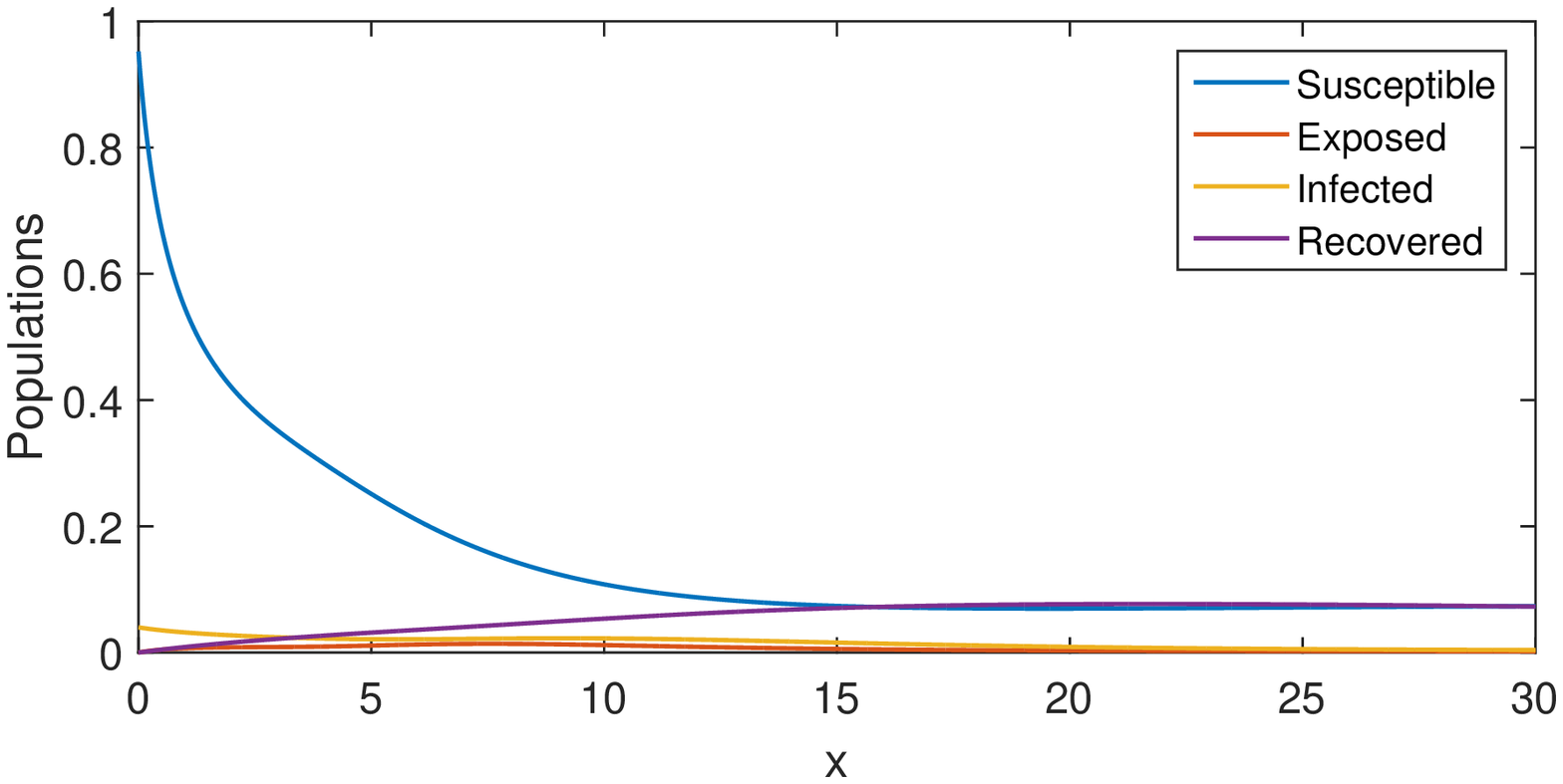}}~
	\subfigure[]{\includegraphics[height=1.7in ,width=2in]{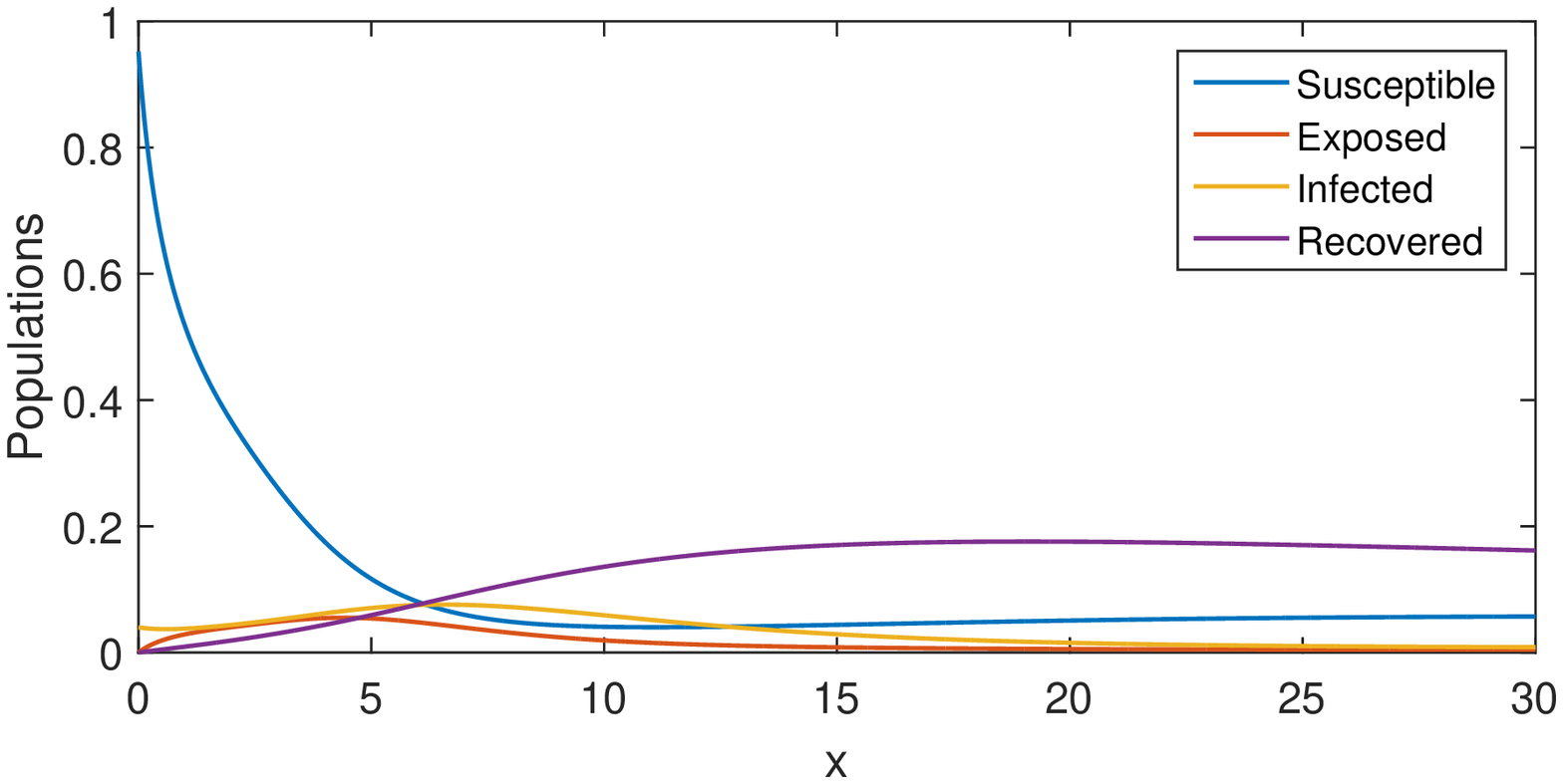}}~
	\subfigure[]{\includegraphics[height=1.7in ,width=2in]{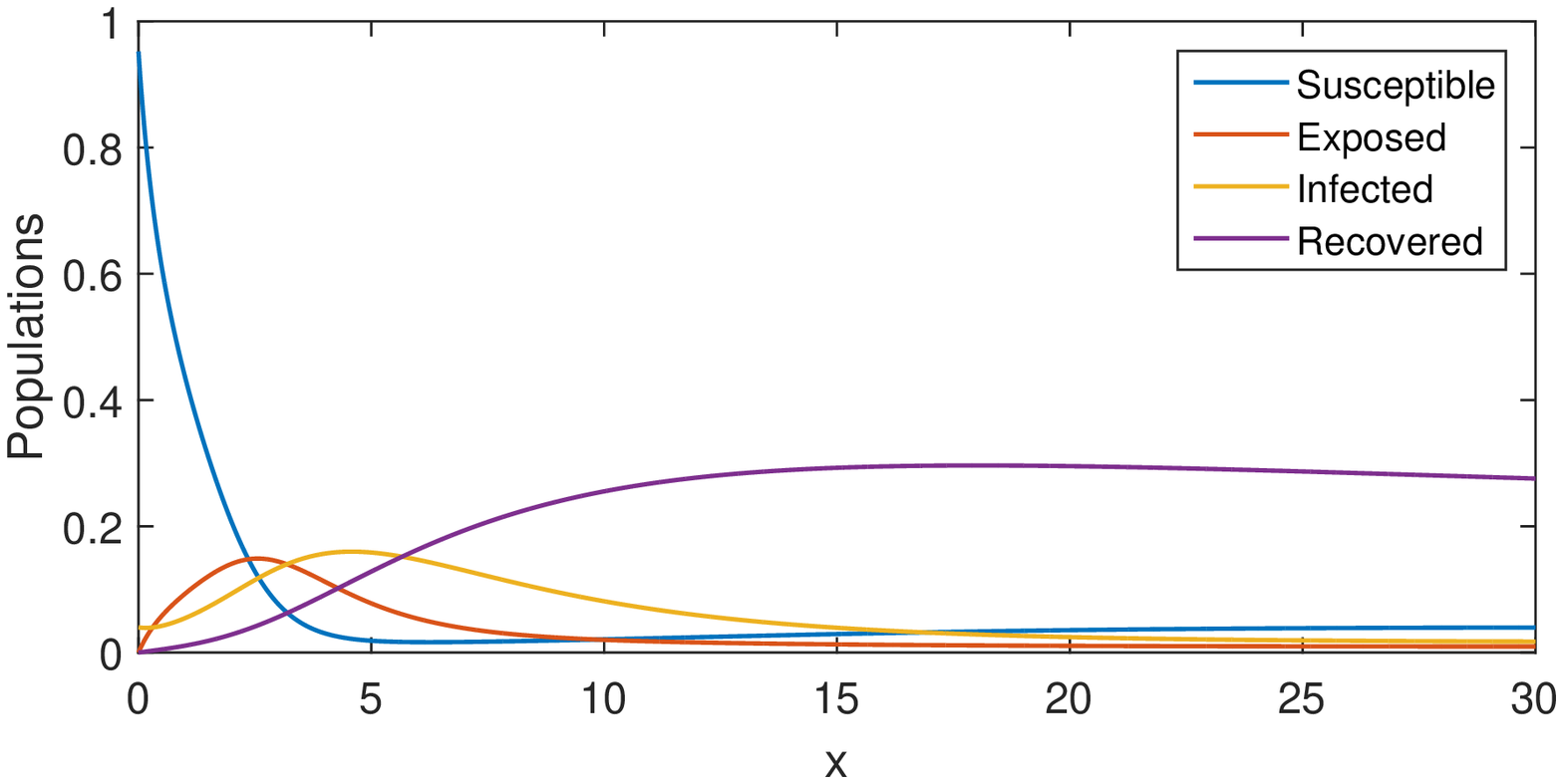}}
	\caption{Time variation of the obtained numerical solutions from (AP)-scheme with $\varepsilon=10^{-6}$ using initial condition $ii)$ and with diffusion, at $x = 0$, for the transmission rate values $\beta=0.03,\,0.075,\,1.12,\,0.1799,\,0.7497,\,2.2491$, the corresponding reproduction ratio is $R_0 = 0.2,\,0.5,\,0.8,\,1.2,\,5,\,15$  }
	\label{F2}
\end{figure}

To well demonstrate the sensitivity of the evolution of the individuals with respect to the transmission rate, we consider different constants values $\beta=0.03,\,0.075,\,1.12,\,0.1799,\,0.7497,\,2.2491$ (the corresponding reproduction ratio is $R_0 = 0.2,\,0.5,\,0.8,\,1.2,\,5,\,15$).  In Figure \ref{F2}, we show the time variation of susceptible, exposed and infected individuals for different values of the transmission rate $\beta$ with diffusion at $x=0$ using the initial condition $i)$. We notice that for small values of the transmission rate, the proportion of the infected population is small. The steady-state ends up with a relatively small proportion of the population in the compartment $R$, while the main proportion of the population remains in the susceptible compartment $S$ (did not catch the disease). However, for relatively higher and moderate values of $\beta$, an important proportion of the population ends up (at the steady-state) in the compartment $R$ (i.e. most of the individuals of the population caught the disease and have been infected then they have recovered). In this case, just a relatively small proportion of the population remains in the compartment $S$. We notice also that the infected and exposed individuals vanish after a reasonable amount of time, while the susceptible and the recovered individuals reach a non zero constant steady-state value. \\

\begin{figure}[h!]
	\centering
	\subfigure[]{\includegraphics[height=1.7in ,width=2.5in]{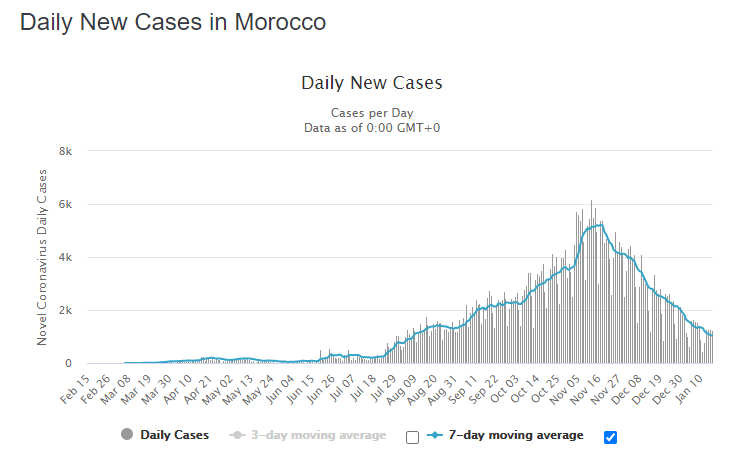}}~
	\subfigure[]{\includegraphics[height=1.7in ,width=2.5in]{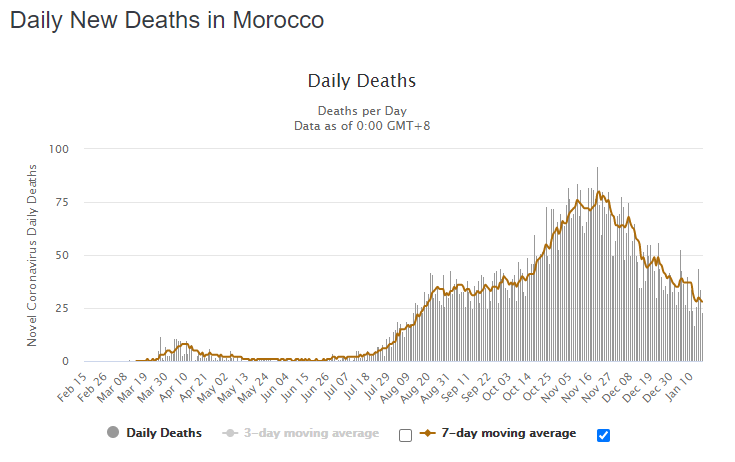}}~
	%\subfigure[]{\includegraphics[height=1.7in ,width=2in]{R04.eps}}~
	%\subfigure[]{\includegraphics[height=1.7in ,width=2in]{R05.eps}}~
	%\subfigure[]{\includegraphics[height=1.7in ,width=2in]{R06.eps}}
	\caption{Daily new cases and deaths in Morocco up to November 17, 2020 \cite{[Rf2]}}
	\label{F5}
\end{figure}

\begin{figure}[h!]
	\centering
	\subfigure[]{\includegraphics[height=1.7in ,width=2in]{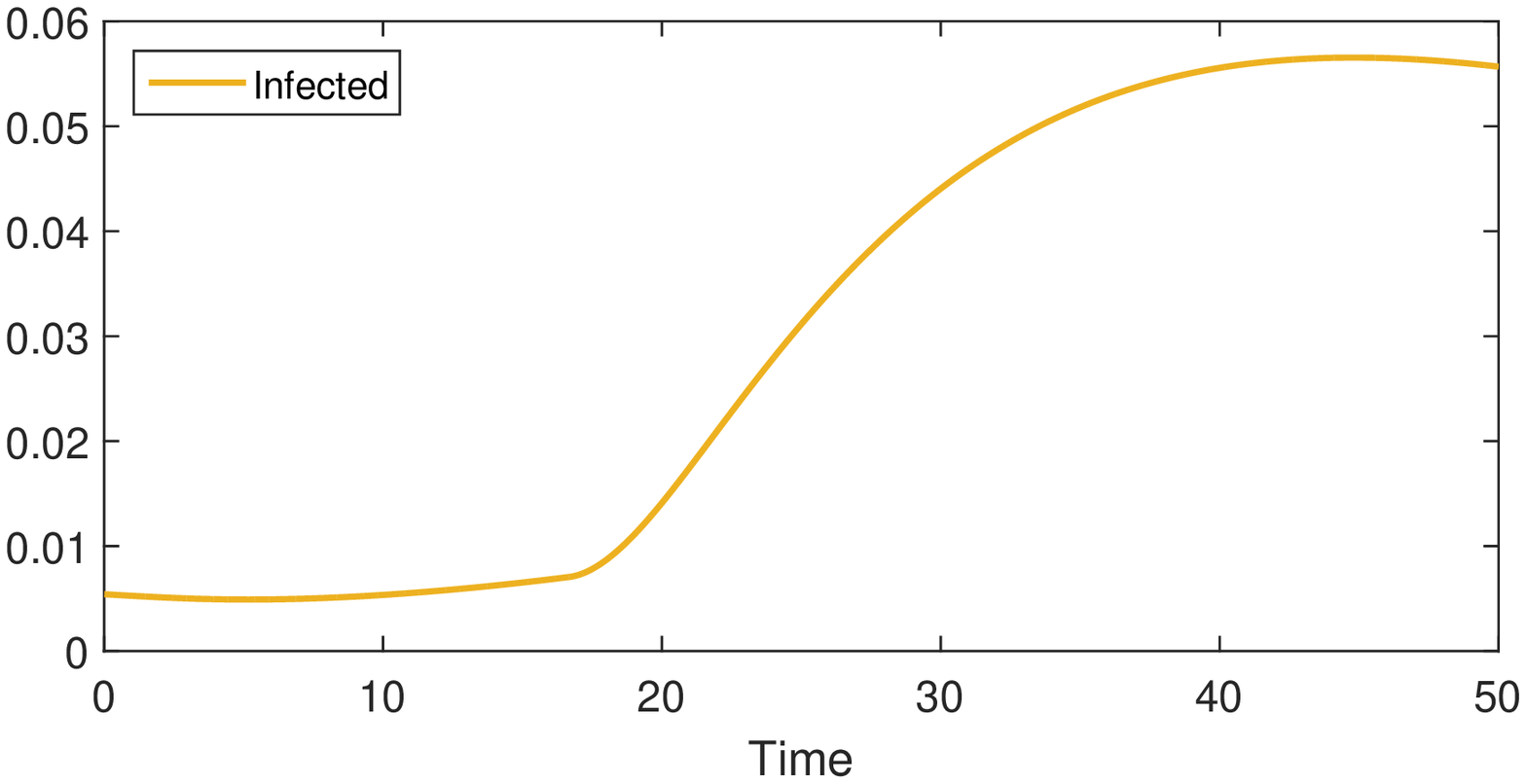}}~
	\subfigure[]{\includegraphics[height=1.7in ,width=2in]{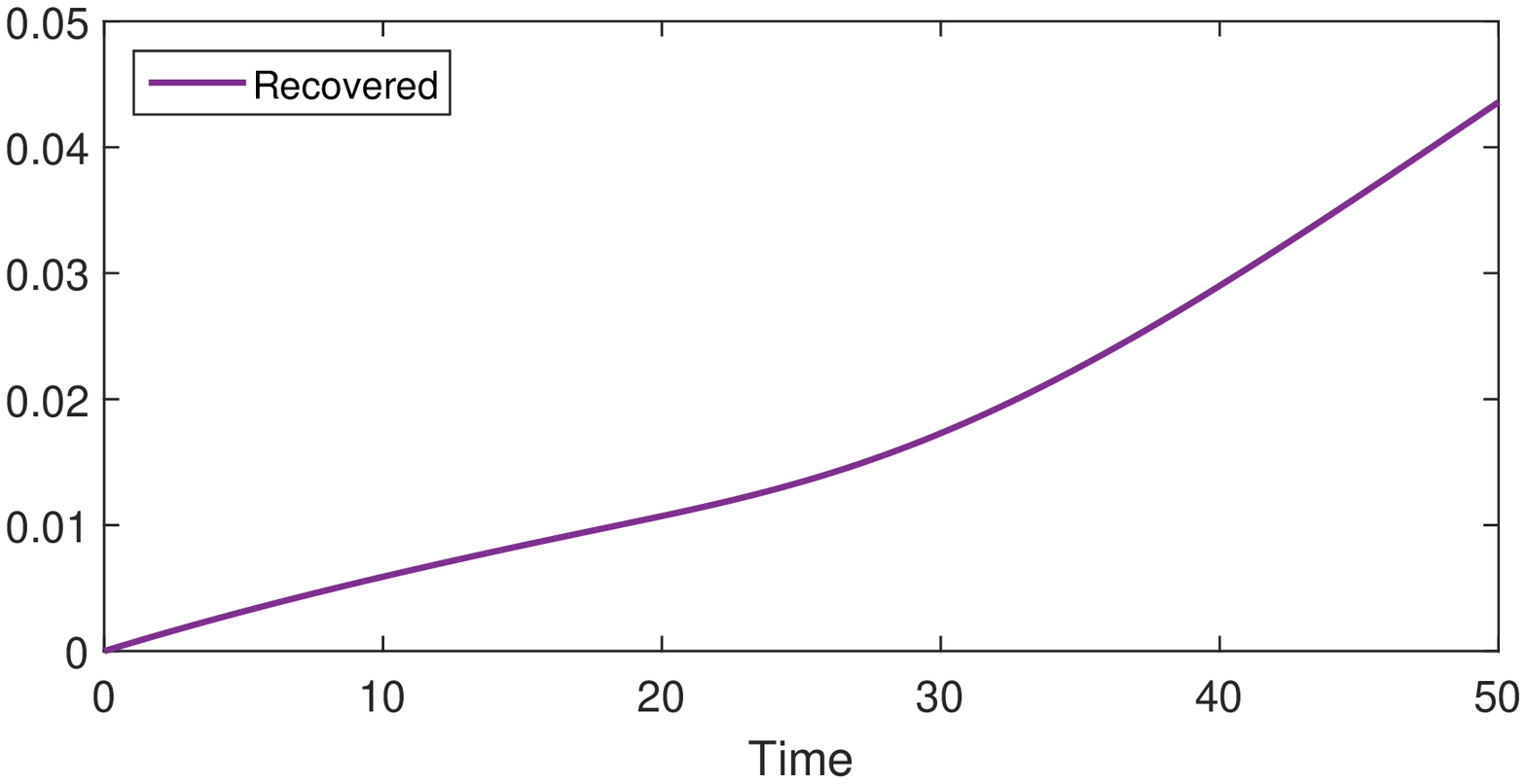}}~
	\subfigure[]{\includegraphics[height=1.7in ,width=2in]{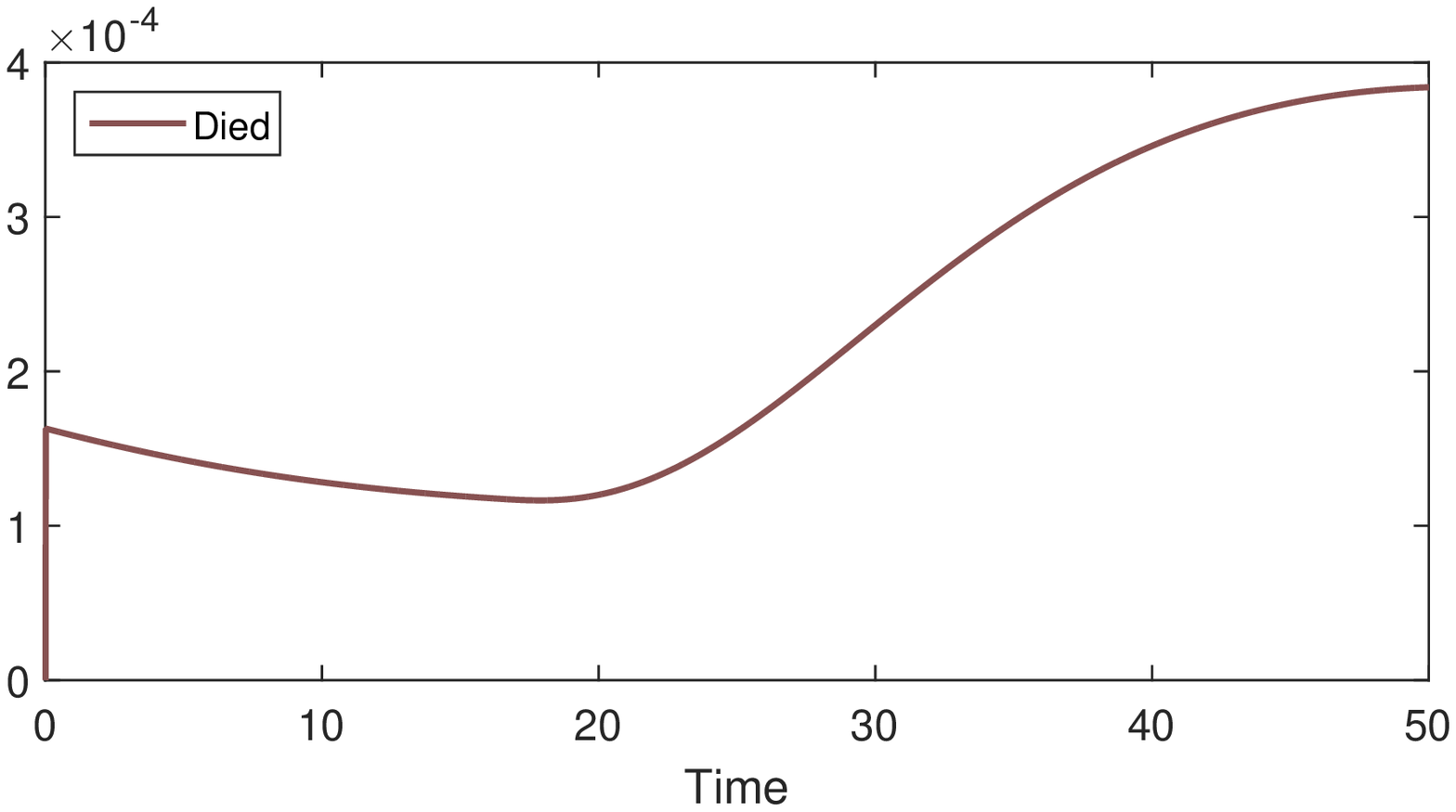}}
	%\subfigure[]{\includegraphics[height=1.7in ,width=2in]{R04.eps}}~
	%\subfigure[]{\includegraphics[height=1.7in ,width=2in]{R05.eps}}~
	%\subfigure[]{\includegraphics[height=1.7in ,width=2in]{R06.eps}}
	\caption{Time variation of the populations of infected, recovered and died individuals obtained from (AP)-scheme with $\varepsilon=10^{-6}$ using initial condition $i)$ with diffusion, at $x = 0.5$, for the transmission rate $\beta(t)$ given by Eq. \eqref{bita1}}
	\label{F4}
\end{figure}

\begin{figure}[h!]
	\centering
	\subfigure[]{\includegraphics[height=1.7in ,width=2in]{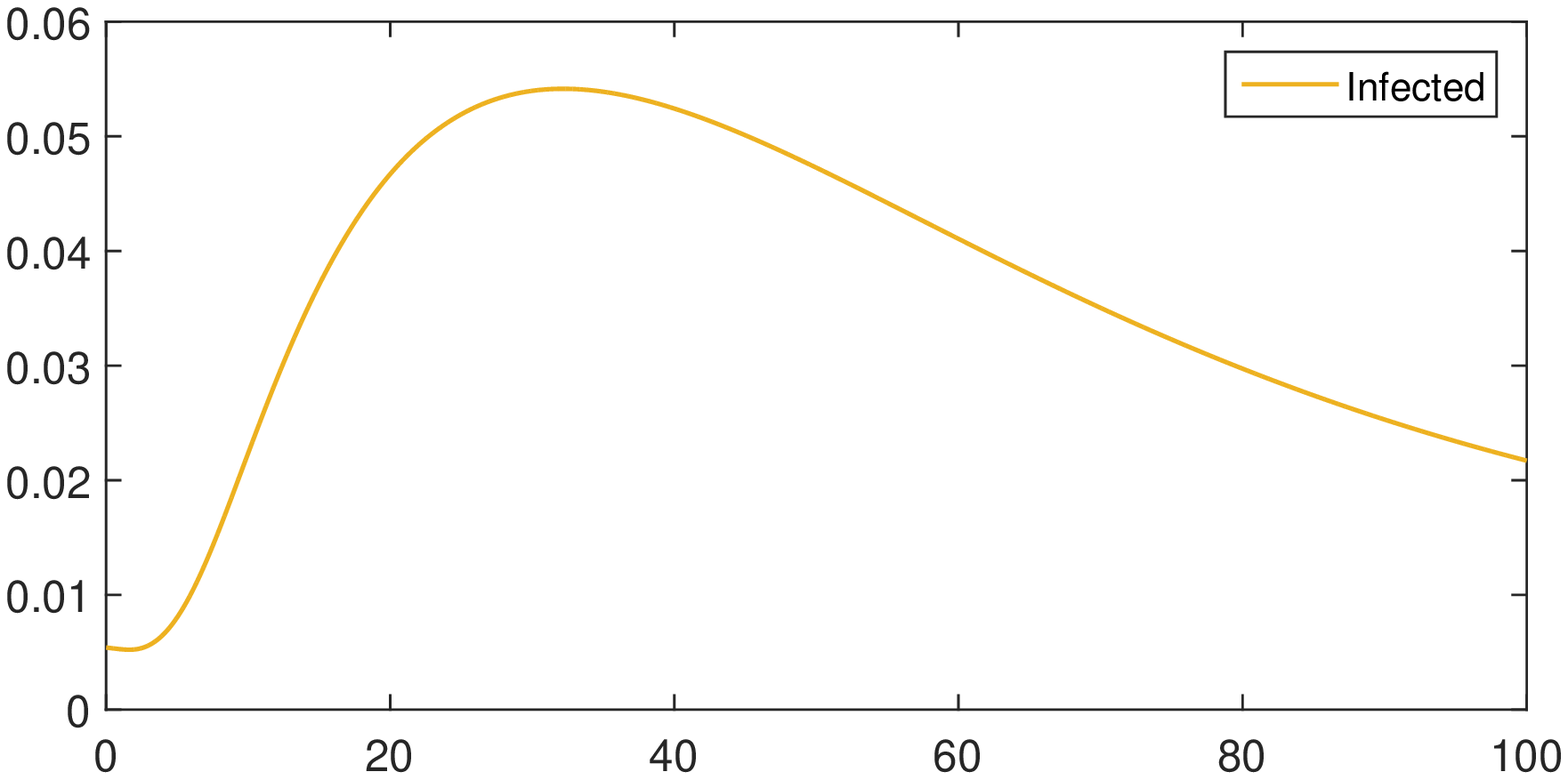}}~
	\subfigure[]{\includegraphics[height=1.7in ,width=2in]{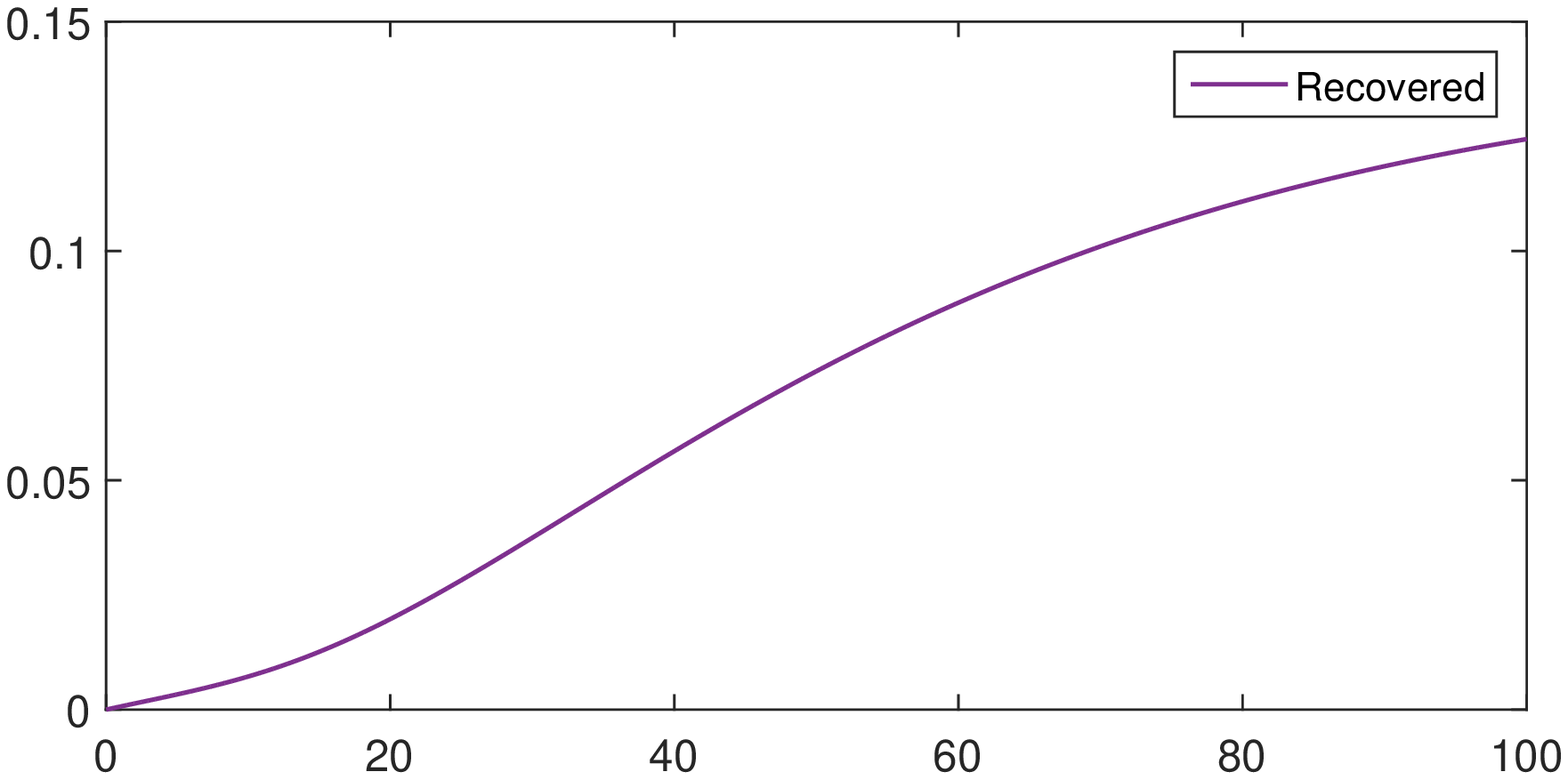}}~
	\subfigure[]{\includegraphics[height=1.7in ,width=2in]{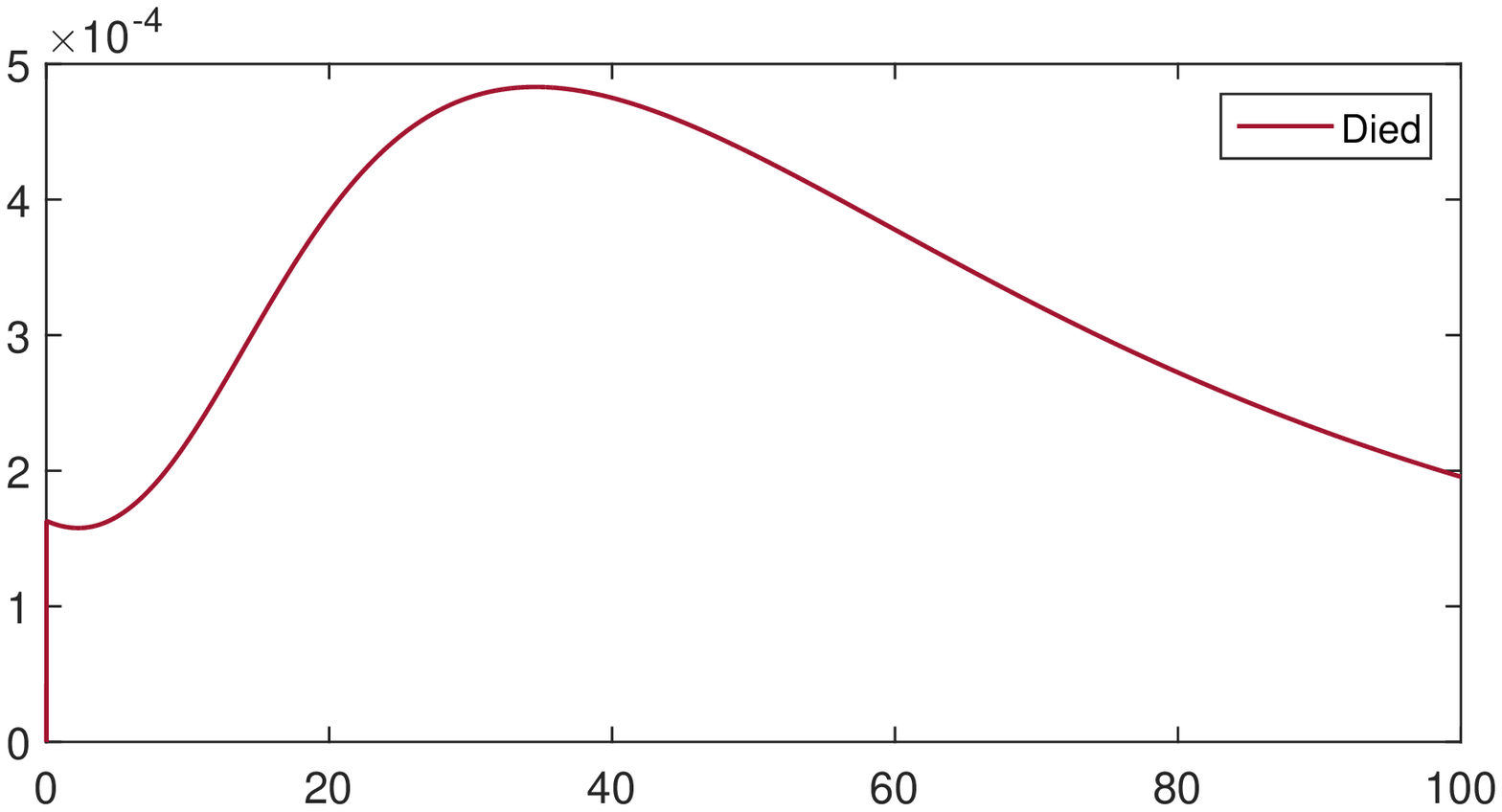}}
	%\subfigure[]{\includegraphics[height=1.7in ,width=2in]{R04.eps}}~
	%\subfigure[]{\includegraphics[height=1.7in ,width=2in]{R05.eps}}~
	%\subfigure[]{\includegraphics[height=1.7in ,width=2in]{R06.eps}}
	\caption{Time variation of the populations of infected, recovered and died individuals obtained from (AP)-scheme with $\varepsilon=10^{-6}$ using initial condition $i)$ with diffusion, at $x = 0.5$, for the transmission rate $\beta(t)$ given by Eq. \eqref{bita2}}
	\label{F6}
\end{figure}

In Morocco, on August 01, 2020 the number of active cases was 25,015 cases and 367 deaths, but as of January 17, 2021 this number increases so fast (458,865 infected and 7,911 deaths). It can be classified into two phases, namely before and after Eid Aldha which was on July 31, 2020. In this religious celebration, peoples used to spend it with their families and visited neighbours. So, many people traveling from city to other ones. On the other hand, unfortunately people do not respect the authority's orders, for instance the partial lockdown in some cities, the obligation to wear a mask and social distance. This leads to a high number of infected and deaths, see Figure \ref{F5}. This situation can be modeled by considering a time-dependent transmission rate $\beta(t)$. We assume the following step-wise functions
\begin{equation}\label{bita1}
\beta(t)=0.075\mathbb{1}_{[0,T/2]}(t)+1.4995\mathbb{1}_{]T/2,T]}(t), 
\end{equation}   
where $T=50$, and
\begin{equation}\label{bita2}
\beta(t)=0.075\mathbb{1}_{[0,T/3]}(t)+1.4995\mathbb{1}_{]T/3,2T/3[}(t)+0.05\mathbb{1}_{]2T/3,T]}(t), 
\end{equation}   
where $T=100$.
It is very clear that the variation in the curves of the infected and deceased populations in Figure \ref{F4} corresponds to the situation in Morocco described above. In Figure \ref{F6}, we notice that the number of infected and deceased populations decreases, while the number of people recovered increases at time $ T>50 $. This happened thanks to the choice of the transmission rate function $ \beta (t) $ given by Eq. \eqref{bita2} where we take a small value of $ R_0 = 0.2 <1 $. Then, we believe that the partial lockdown in Moroccan applied from December 24,  2020 and the national vaccination campaign against COVID-19 in the upcoming weeks stated from the Royal Office announced November 9, 2020 will gradually eliminate the virus over the next few months.
%\begin{figure}[h!]
%	\centering
%	\subfigure[]{\includegraphics[height=1.7in ,width=2in]{ST1.eps}}~
%	\subfigure[]{\includegraphics[height=1.7in ,width=2in]{ET1.eps}}~
%	\subfigure[]{\includegraphics[height=1.7in ,width=2in]{IT1.eps}}
%		\subfigure[]{\includegraphics[height=1.7in ,width=2in]{ST2.eps}}~
%	\subfigure[]{\includegraphics[height=1.7in ,width=2in]{ET2.eps}}~
%	\subfigure[]{\includegraphics[height=1.7in ,width=2in]{IT2.eps}}
%		\subfigure[]{\includegraphics[height=1.7in ,width=2in]{ST3.eps}}~
%	\subfigure[]{\includegraphics[height=1.7in ,width=2in]{ET3.eps}}~
%	\subfigure[]{\includegraphics[height=1.7in ,width=2in]{IT3.eps}}
%		\subfigure[]{\includegraphics[height=1.7in ,width=2in]{ST4.eps}}~
%	\subfigure[]{\includegraphics[height=1.7in ,width=2in]{ET4.eps}}~
%	\subfigure[]{\includegraphics[height=1.7in ,width=2in]{IT4.eps}}
%	\caption{Dynamics of the densities $S$ (first column), $E$ (second column) and $I$ (third column) obtained from the (AP)-scheme with $\varepsilon =2\times10^{-k}$, $k = 0,\,1,\, 2,\,3,\, 4,\, 6$ against the macroscopic scheme with initial conditions $i)$ at successive time $t=0.5,\,1,\,5,\,10$.}
%	\label{F0}
%\end{figure}

%\subsection{Parameter estimation, model validation and prediction}


\section{Conclusion and perspectives}

In this paper, a time-independent SEIRD reaction-diffusion system for individual populations has been proposed. This system has been derived from kinetic model (\ref{mM1}) by using the micro-macro decomposition method. It has been shown that the proposed (AP)-scheme is uniformly stable along the transition from kinetic to macroscopic regimes. Various promising numerical simulations have been provided. Specifically, it has shown that the presence of the diffusion terms in system \eqref{Cross-Diffusion} influences the spreading of the pandemic. On the other hand, the sensitivity to the transmission rate is demonstrated. Indeed, for small values of the transmission rate $\beta$, the proportion of the infected population is small, and the steady-state ends up with a relatively small proportion of the population in the compartment $R$, while the main proportion of the population remains in the susceptible compartment $S$. While, for relatively moderate and higher values of $\beta$, an important proportion of the population ends up at the steady-state in the compartment $R$, and the infected and exposed individuals vanish after a reasonable amount of time, while the susceptible and the recovered individuals reach a non zero constant steady-state value. 
Finally, the importance in considering a time-dependent rate transmission has been demonstrated. The obtained numerical results describe the Moroccan actual situation where the go back to partial lockdown decreases the high numbers of the infected cases and deaths. However, the vaccination campaign, which will be implemented soon, will gradually eliminate the virus over the next few months. We believe that this paper opens such interesting perspectives: Numerical study of the proposed system for the real empirical data, extension of the macroscopic model by considering a time-space diffusions $d_i(t,x)$ and the rate transmission $\beta(t,x)$. Moreover, we think that it is interesting to apply the presented method in this paper to other applications, for instance mechanism of normal and infected cells. An interesting macroscopic models can be found in the papers by \cite{[BPTW19],[SAMH13]}. \\

%{\footnotesize
	%\nocite{*}
%	\bibliographystyle{elsarticle-num}
%	\bibliography{BIB}}

\end{document}